%% file: main.tex
\documentclass[english,11pt]{elsarticle} 
\usepackage[T1]{fontenc}
\usepackage{float}
\usepackage{amsmath}
\usepackage{amsthm}
\usepackage{amssymb}
\usepackage{amsfonts}
\usepackage{bbm}
\allowdisplaybreaks
\usepackage{graphicx}
\usepackage{natbib}
\setcitestyle{square, comma, numbers,sort&compress, super}
\usepackage{float}
\usepackage{subfloat}
\usepackage{epstopdf}
\usepackage{hyperref}
\hypersetup{colorlinks,allcolors=black}
\usepackage{nomencl}
\makenomenclature
\usepackage{algorithm}
\usepackage{algorithmic}

\usepackage{url}
\usepackage{comment}
\usepackage{ragged2e}
\usepackage{pgfplots}
\pgfplotsset{compat=newest}
\usetikzlibrary{plotmarks}
\usetikzlibrary{arrows.meta}
\usepgfplotslibrary{patchplots}
\usepackage{grffile}
\usepackage{amsmath}
\usepackage{autobreak}
\usepackage{appendix} 
\usepackage{appendix}
\usepackage[symbol]{footmisc}

\usepackage{paralist}
\usepackage{tablefootnote}

\input{format}

\biboptions{numbers,sort&compress} 

\begin{document}
\begin{frontmatter}
 \title{\textbf{Regulating Transportation Network Companies with a Mixture of Autonomous Vehicles and For-Hire Human Drivers}}	

\author{Di Ao}
\ead{daoaa@connect.ust.hk}
\author{Jing Gao}
\ead{jgaoax@connect.ust.hk}
\author{Zhijie Lai}
\ead{zlaiaa@connect.ust.hk}
\author{Sen Li\corref{mycorrespondingauthor}}
\cortext[mycorrespondingauthor]{Corresponding author}
\ead{cesli@ust.hk}
\address{Department of Civil and Environmental Engineering, The Hong Kong University of Science and Technology}

\begin{abstract}
This paper investigates the equity impacts of autonomous vehicles (AV) on for-hire human drivers and passengers in a ride-hailing market, and examines regulation policies that protect human drivers and improve transport equity for ride-hailing passengers. We consider a transportation network companies (TNC) that employs a mixture of AVs and human drivers to provide ride-hailing services. The TNC platform determines the spatial prices, fleet size, human driver payments, and vehicle relocation strategies to maximize its profit, while individual passengers choose between different transport modes to minimize their travel costs. A market equilibrium model is proposed to capture the interactions among passengers, human drivers, AVs, and TNC over the transportation network. The overall problem is formulated as a non-concave program, and an algorithm is developed to derive its approximate solution with a theoretical performance guarantee. Our study shows that TNC prioritizes AV deployment in higher-demand areas to make a higher profit. As AVs flood into these higher-demand areas, they compete with human drivers in the urban core and push them to relocate to suburbs. This leads to reduced earning opportunities for human drivers and increased spatial inequity for passengers. To mitigate these concerns, we consider: (a) a minimum wage for human drivers; and (b) a restrictive pickup policy that prohibits AVs from picking up passengers in higher-demand areas. In the former case, we show that a minimum wage for human drivers will protect them from the negative impact of AVs with negligible impacts on passengers. However, there exists a threshold beyond which the minimum wage will trigger the platform to replace the majority of human drivers with AVs. In the latter case, we show that prohibiting AVs from picking up passengers in higher-demand areas not only improves the spatial equity of ride-hailing services for passengers, but also substantially increases human driver surplus and restricts the increase of total fleet size compared to the unregulated case. These results are validated with realistic case studies for San Francisco.
\end{abstract}

\end{frontmatter}

\section{Introduction}
\label{section1}
Transportation network companies (TNCs), such as Uber, Lyft, and Didi, have deeply disrupted urban transportation \cite{yang2020integrated,li2021impact,wei2020mixed,ozkan2020joint,sun2019optimal,karamanis2020assignment}. They provide mobility-on-demand services that immediately respond to passengers requests from anywhere at any time \cite{oh2020assessing}. These compelling features make TNC a booming market: 
for many travelers, TNC has become an indispensable transport option for commuting and entertaining trips. However, despite their popularity, many TNC platforms are still struggling to find their ways to make profits. The profitability challenge of TNCs is deeply  rooted in its business model: to ensure a short response time, they have to recruit a large number of available but idle drivers. This incurs additional operational costs which have to be shared between the TNC platforms and the for-hire drivers, creating serious labor tensions in the ride-hailing industry. 
As reported in \cite{parrott2020minimum}, the net payment for drivers in Seattle is about \$\,9.73 after excluding the expenses on vehicle registration, operating costs, fuel costs, insurance costs, etc, which is far less than the \$\,16.39 hourly minimum wage of the city. Another work \cite{parrott2018earnings} investigated the working conditions of for-hire drivers in New York City and found that a majority of for-hire drivers in NYC works full-time and that 85 percent of the drivers make less than the city's minimum wage. The concerns on the working conditions of for-hire drivers have prompted multiple cities (e.g., Seattle, NYC) to implement minimum wage standards on the ride-hailing industry, which makes it even harder for the platforms to find solutions towards profitability.

The advance of autonomous-driving technology brings significant potential for the TNC platform to reduce its reliance on for-hire human drivers and improve its financial viability \cite{wang2020control}. Although Level-5 self-driving technology is still under development, Level-4 AVs operated within geofenced areas under favorable weather conditions have already been deployed in various cities nowadays (e.g., Phoenix, AZ) \cite{Waymo}. These Level-4 AVs can be jointly managed with a group of for-hire human drivers, comprising a mixed fleet, to offer ubiquitous mobility-on-demand services to TNC passengers from all locations at all times in the near future. Compared to human drivers, AV-enabled mobility services have several advantages for the TNC platform. On the one hand, AV fleet will remove the cost associated with human drivers, which accounts for the majority of the platform's operational cost. This can significantly improve the profitability of the TNC platform. On the other hand, AVs can be centrally dispatched by the platform to relocate to preferable locations. This enables the platform to mitigate the spatial mismatch between supply and demand with higher efficiency at a lower cost. Both changes will significantly benefit TNC platforms, enabling them to provide higher-quality services while enjoying higher profits.

While the deployment of AVs benefits TNC platforms, it may also open a Pandora's box of negative externalities, leading to unfavorable market outcomes for human drivers and TNC passengers. For drivers, a prevalent concern is that AVs will compete with human drivers and reduce their hourly wage, or even worse, they may completely replace human drivers in the TNC market and lead to massive loss of job opportunities \cite{hancock2019future}.  {This may provoke social instability because more than 80\% of full-time drivers purchase their vehicles just to enter the TNC industry \cite{parrott2018earnings}.} For passengers, the prevalent concern is that TNC platforms tend to concentrate their services in highest-density areas of the city. This is the natural outcome of a free market as TNC platforms often prioritize the deployment of their resources in more profitable areas when they seek to maximize their profits. However, these market incentives not only contribute to the severe traffic congestion in urban core, but also incur inadequate services in suburbs and city outskirts. Unfortunately, this inequitable distribution of TNC services will be further exacerbated as the proliferation of autonomous-driving technology in the ride-hailing industry enables the TNC platform to enjoy lower cost and better control of their fleet. Without proper regulation, market incentives of the TNC platforms may drive us to a disastrous mobility future, where for-hire human-drivers are massively replaced with AVs, transportation-disadvantaged groups in underserved areas are further left behind, and the benefits of AVs are disproportionately allocated to TNC platforms instead of those who need them most. These concerns motivate us to examine the future TNC market with an equity focus, raising the following research questions:
\begin{inparaenum}[(i)]
\item How will AVs affect human drivers over the transportation network?
\item What will AVs affect the spatial equity of TNC services for passengers?
\item How to protect human drivers and passengers from the negative impacts of AVs through regulations?
\end{inparaenum}

This paper investigates the aforementioned questions with a mathematical model that captures the interactions between AVs, human drivers, passenger, and the TNC platform in a regulated TNC market. In particular, we consider a TNC market with a platform that recruits a mixture of AVs and human drivers to offer mobility-on-demand services. In this market, passengers choose between TNC and other transport modes to minimize their travel costs, human drivers make market entry and vehicle relocation decisions to maximize their earnings, while the TNC platform determines the spatial prices, the fleet size, and AV relocation strategies to maximize its profit. An economic equilibrium\footnote{The ``equilibrium" refers to market outcomes arising from the interactions between passengers, AVs, human drivers, and the TNC platform.} model is formulated to characterize the incentives of passengers, AVs, human drivers, and the TNC platform over the transportation network. At the equilibrium solution, the spatial distributions of passenger and drivers are derived, the network flows of AVs and human drivers are balanced, and the platform profit is maximized. Our major contributions are summarized below:
\begin{enumerate}[(1)]
    \item We propose an economic analysis for a regulated ride-hailing network with mixed autonomy, which not only captures the incentives of passengers, AVs, human drivers, the platform, and the policy maker, but also characterizes their interactions over the transportation network. We prove that the solution to the proposed network equilibrium model exists, and we developed a decomposition-based algorithm to compute the approximate solution to the platform's profit maximization problem with a theoretical performance guarantee. In particular, an upper bound is derived to characterize the optimality gap between the derived solution and the globally optimal solution to the platform's non-concave profit maximization problem. 
    \item We evaluate the impacts of AVs on human drivers and passengers in the ride-hailing network. We show that the TNC platform prioritizes AV deployment in the urban core as there is higher passenger demand. This creates competition between AVs and human drivers in these areas and pushes human drivers to relocate to suburbs, which significantly reduces the earnings and job opportunities of human drivers. Furthermore, the geographical concentration of TNC services enlarges the existing inequity gap among TNC passengers in different areas of the city. This indicates that the benefits of AVs will be unfairly distributed among passengers, and regulatory policies are needed to improve transport equity for disadvantaged groups in underserved areas.
    \item We investigate the impacts of a wage floor on human-driver earnings in the mixed environment. We show that a carefully selected minimum wage on human driver earnings can benefit drivers by increasing the driver wage and creating more job opportunities. In the meanwhile, it does not incur significant loss to TNC passengers. However, there exists a threshold beyond which the minimum wage will trigger a paradigm shift of labor supply and the platform will replace the majority of human drivers with AVs. This indicates that minimum wage for human drivers should be exercised with discretion in the mixed environment, because an overly protective wage floor can lead to massive loss of job opportunities. 
    \item { We evaluate the impacts of a restrictive pickup policy that prohibits AVs from picking up passengers in higher-demand areas (they can drop off passengers in these areas), which is consistent with existing policies in New York City and Hong Kong, where taxis are classified by colors and certain color of taxis can not pick up passengers in urban centers. We show that prohibiting AVs from picking up passengers in higher-demand areas not only improves the
    spatial equity of TNC services for passengers, but also substantially increases the surplus of human driver and restricts the increase of total fleet size compared to the unregulated case. However, although the proposed policy benefits passenger in lower-demand areas, it will negatively affect passenger in higher-demand areas, which reflects a trade-off between equity and service quality.}
\end{enumerate}

\section{Related Works}
This section explores related works in the existing literature and discusses how our paper differs from them. We will cover the management of ride-sourcing platforms, the management of mixed fleets, and the regulation on ride-sourcing platforms and AMoD services, respectively.  
\label{Section2}
\subsection{Management of TNC Platforms}
The pricing strategy is an important research topic in the study of TNC platforms as it is crucial for balancing the supply and demand in a ride-hailing network \cite{zhu2021mean,chen2020dynamic,yang2020optimizing,candogan2012optimal,bimpikis2019supply,ke2017short,zhang2021pool}. This is one of the key factors to differentiate between the ride-sourcing industry and the traditional taxicab industry whose ride fare is strictly regulated by the government. One stream of works studied spatial pricing, where the platform can determine the location-based ride fares and driver payments to facilitate better vehicle repositioning. For instance,  in \cite{bimpikis2019spatial}, spatial price discrimination is considered in the context of  a ride-sourcing platform that serves a network of locations. It revealed that profits and consumer surplus at the platform’s optimal pricing and compensation policy are maximized when the demand pattern is balanced across the transportation networks. In \cite{zha2018geometric}, the effects of spatial pricing on ride-sourcing markets are evaluated based on a discrete time geometric matching framework. It was reported that the platform may resort to relatively higher price to avoid the inefficient supply of vehicles if spatial price differentiation is not permitted. In \cite{chen2019optimal}, the authors investigate a long-term optimal spatial pricing strategy for ride-sourcing platform and shows that the optimal number of participating drivers and their optimal wages will be influenced not only by the pricing strategy but also by the levels of service of the service zones. Aside from spatial pricing, another stream of works considers surge pricing, where the platform dynamically adjust the prices according to real-time imbalances between supply and demand to improve customer experience and increase platform profit. For instance, in \cite{guda2019your}, the authors considered the use of surge pricing in ride-sourcing platforms and showed that surge pricing can be useful even in zones where supply exceeds demand. In \cite{ma2020spatio}, a spatiotemporal pricing mechanism was proposed for the ride-sourcing platform, which is welfare-optimal, envy-free, individually rational, budget balanced and core-selecting in equilibrium. In \cite{afifah2022spatial}, the authors investigated the issue of spatial pricing in the context of ride-sourcing services. They found that the problem of determining optimal pricing can be effectively addressed through a convex reformulation, provided that the waiting time for passengers was relatively small compared to the travel time. {\em We emphasize that all aforementioned works differ from our paper as they primarily focus on human drivers, which neglect the interactions between AVs and human drivers in the TNC market.} 


\subsection{Management of Mixed Fleet}
TNC platforms can also deploy AVs or a mixture of AVs and human drivers to offer autonomous mobility-on-demand (AMoD) services \cite{pavone2015autonomous,zhang2015models,zhang2016control, wallar2019optimizing}. In \cite{zhang2016control}, a queuing-theoretical method was proposed for the modeling, analysis, and control of AMoD systems, and an optimal rebalancing algorithm were developed to minimize the number of rebalancing vehicles while keeping vehicle availability balanced throughout the network. In \cite{wallar2019optimizing}, an offline method was presented to optimize the vehicle distribution and fleet size based on historical demand data for AMoD systems under ride-sharing. An algorithm was developed to determine how many vehicles are needed, where they should be initialized, and how they should be routed to serve the travel demand. In the case of mixed traffic, many works have investigated the interactions between AVs and human-driven vehicles through cooperative adaptive cruise control system or in signalized intersection \cite{qin2021lwr,jin2020analysis,liu2021latency}. Some other works focus on the optimal pricing and control of for-hire human drivers and AVs on the ride-hailing network. In \cite{lokhandwala2018dynamic}, agent-based simulation was employed to study a TNC platform with a mixed fleet, which includes traditional taxis and shared autonomous taxis, to analyze the advantages and disadvantages of taxi sharing networks. In \cite{wei2019ride,wei2020mixed}, the interaction between human-driven vehicles and autonomous vehicles under a multi-location equidistant network is researched to maximize the platform profit. It showed that in some cases, there is a regime in which the platform chooses to mix AVs with human drivers to maximize its profit, while in other cases, the platform will use only human drivers or AVs depending on the relative cost of AVs. In \cite{yang2020planning,yang2021real}, an AMoD system with a mixed fleet was modeled as a Stackelberg game where the platform serves as the leader and human-driven vehicles serve as the followers. A steady-state model was proposed to determine the planning variables, and a time-varying model was formulated to design a real-time coordination algorithm for AVs. In \cite{mo2022modeling}, a monopoly TNC market with mixed fleets of AV and human drivers are considered, and the impacts of congestion externality are explicitly characterized in an aggregated market equilibrium model. In \cite{xie2022two}, a multi-agent reinforcement learning framework is developed for the dynamic management of the mobility-on-demand platform with a mixed fleet of AVs and human drivers. The algorithm addresses the interaction between the platform and self-interested human drivers through a two-sided design, and the proposed framework is demonstrated in a case study for New York City. 

{\em Our paper differs from all aforementioned works as we focus on equity analysis and regulatory policies.  Different from existing works, which mainly studies the benefits of AVs to the profit-maximizing platform, we quantify its negative impacts on for-hire human drivers and identify its potential damage to transport equity. Furthermore,  we evaluated the impacts of various regulations that protects human drivers and improve transport equity for TNC passengers in the mixed environment. These elements are crucial, as they not only ensure that the benefits of AVs reach disadvantaged groups in underserved areas, but also generates important insights for the government to make forward-looking policies. However, these elements are missing in all existing literature.}

\subsection{Regulations on TNC Platforms}
The explosive growth of TNCs in recent years imposes significant negative externalities to  road users and gig workers, which inspires an emerging stream of works that focuses on TNC regulations. In \cite{parrott2018earnings,parrott2020minimum}, a minimum wage floor standard was proposed for TNC drivers in NYC and Seattle. Data were collected to estimate the net earning of for-hire drivers in these cities, and possible consequences of the proposed minimum wage standard was investigated. In \cite{li2019regulating}, an economic equilibrium model was proposed to investigate the impacts of various policies on the TNC market, including a driver wage floor, an upper bound on the number of drivers, and a congestion charge on TNC trips. It showed that the minimum wage for TNC drivers may benefit both drivers and passengers and promote the efficiency of the entire system at the cost of TNC profit. Reference \cite{li2021impact} formulated an analytical model to evaluate the joint impacts of congestion charge and minimum wage for TNC drivers, and showed that the minimum wage may mitigate the impact of congestion charges on the number of TNC vehicles. In \cite{li2021spatial}, a spatial equilibrium model was proposed to investigate the impact of cordon-based congestion charge on the ride-hailing network. In \cite{zhang2021mitigating}, the authors considered a TNC platform that offers both solo and pooling e-hail services and investigated the impacts of various congestion mitigation policies, such as a trip-based charge, a cordon-based charges, and a cruising cap that requires the TNC to maintain the fleet utilization ratio in the city center above a threshold. Other regulation policies, including limiting the total number of vehicles, capping the platform commission rate, etc., were discussed in \cite{yu2020balancing,zha2016economic,vignon2021regulating}.

Compared to TNCs, regulations on AMoD systems are relatively under-examined. In \cite{freemark2020policies},  future policies on AV regulations were surveyed among multiple municipal officials from different cities in the U.S.  In \cite{simoni2019congestion}, congestion pricing and road tolling for AMoD services were investigated to moderate the travel demand,  promote social welfare, and incentivize environmentally friendly travel modes. In \cite{salazar2019intermodal}, an intermodal AMoD system is considered, where the AV fleet is operated in combination with public transport. The pricing and tolling mechanisms were designed to achieve social optimum in a perfect market with selfish agents, and numerical studies were presented to demonstrate the benefits of the integrated AMoD-transit model. In \cite{dandl2021regulating}, a tri-level model was formulated to capture the interactions among public-sector decision makers, mobility service providers,  and travelers. Bayesian optimization method was adopted  to solve the optimization problem, and a case study for Munich was presented to validate the proposed methods.

{\em Note that our paper differs from all aforementioned works as we study regulations on TNC markets in the mixed environment, where AVs interact with for-hire human drivers to offer TNC services to TNC passengers. This is an important case because in the foreseeable future, it is very likely that Level-4 AVs will be deployed for commercial use first, well before Level-5 AV becomes available. In this case, the AV fleet has to be complemented with for-hire human drivers to ensure that TNC services are accessible to all passengers, at all times, and at all locations of the city. }

\section{Mathematical Model for the Ride-Hailing Network} 
\label{lowerlevel}

This section formulates a spatial equilibrium model to characterize the ride-hailing network with mixed autonomy under regulations.  The following subsections presents the details of this model.

\subsection{Problem Setup}
Consider a city comprised of $M$ zones. Let $\mathcal{V}$ denote the set of zones in the road network. The TNC platform hires a mixture of AVs and human drivers to serve passengers in the city. Each passenger sends a ride request from the origin $i\in\mathcal{V}$ to the destination $j\in\mathcal{V}$. Upon receiving the request, the platform will dispatch the closest idle vehicle in the same zone to pick up the passenger, no matter it is an AV or a human-driven vehicle. We assume that AVs and human drivers provide the same quality of service and passengers are indifferent to the type of vehicle that serves the trip. After the passenger is picked up, he/she will be transported from $i$ to $j$ following the shortest path, which takes an average time of $t_{ij}$. After the passenger is dropped off, the vehicle will either stay in zone $j$ or relocate to a different zone to seek the next passenger. Below we characterize the incentives of passengers, drivers, and the platform, respectively. 

\subsection{Incentives of Passengers}
In a multimodal transportation system, each passenger may have a few transport options, such as TNC, taxi, public transit, walking, etc.  For a passenger that travels from zone $i$ to $j$, his/her choice depends on the generalized travel costs of each transport mode. Distinct passengers may have different preferences. For instance, the elderly and disabled passengers may prefer to use ride-hailing passengers driven by human drivers instead of AVs. To characterize the heterogeneous preference, we categorize passengers into two classes: (1) Class 1: those without a specific preference for AV or human drivers; and (2) Class 2: those with a preference for human drivers. For simplicity, we assume that the ride fare and travel speed of AVs and human drivers are both identical, and that the first class of passengers are indifferent to AVs and human drivers in making their travel mode choice. In this case, the cost of a TNC ride from zone $i$ to zone $j$ can be defined as the weighted sum of passenger waiting time and trip fare:
\begin{equation}
\begin{cases}
        c_{ij}^{AH} = \alpha {w^{p,AH}_i} + r_i t_{ij},\\
        c_{ij}^{H} =\alpha {w^{p,H}_i} + r_i t_{ij},
\end{cases}
\label{1}
\end{equation}
where $c_{ij}^{AH}$ and $c_{ij}^H$ represent the generalized travel costs for Class 1 and Class 2 passengers, respectively; $w_{i}^{p,AH}$ and $w_{i}^{p,H}$ represent the average passenger waiting times for Class 1 and Class 2 passengers in zone $i$, respectively; $\alpha$ is the passenger's value of time; $r_i$ denotes the average ride fare per unit time for a trip originating from zone $i$; and $t_{ij}$ denotes the average travel time from zone $i$ to zone $j$.
Note that the passenger waiting time $w_{i}^{p,AH}$ and $w_{i}^{p,H}$ depend on the number of idle vehicles in each zone\footnote{We implicitly assume that passengers from zone $i$ will only be matched to an idle vehicle in the same zone. This is a reasonable assumption when the zone of the city is large and the vehicle supply is sufficient so that each passenger will be picked up by a nearby vehicle with high probability. The case of inter-zone matching is omitted in this work for simplicity.}: if there are more idle vehicles, the average passenger waiting time should be shorter. To capture this relation, we denote $N_i^I$ as the number of idle vehicles in zone $i$, which includes both idle AVs $N_i^{I,A}$ and idle human-driven vehicles $N_i^{I,H}$. Note that Class 1 passenger can be matched to either AV or human drivers, whereas Class 2 passenger can only be matched to human drivers. Therefore,  $w_{i}^{p,AH}$ depends on the total number of idle vehicles $N_i^I=N_i^{I,A}+N_i^{I,H}$, and $w_{i}^{p,H}$ depends on the number of idle human drivers $N_i^{I,H}$. With slight abuse of notation, the average passenger waiting time can be represented as decreasing functions $w_{i}^{p,AH}\left(N_i^I\right)$ and $w_{i}^{p,AH}\left(N_i^I\right)$, respectively. We assume that the passenger waiting time  follows the well-established “square root law” (\cite{li2019regulating} and \cite{arnott1996taxi}):
\begin{equation}
    \begin{cases}
         w_i^{p,AH}=\frac{L}{v_i\sqrt{N_i^{I}}}, \\
         w_i^{p,H}=\frac{L}{v_i\sqrt{N_i^{I,H}}}.
    \end{cases}
    \label{30}
\end{equation}
where $L$ is a model parameter that depends on the area of the city and supply-demand distribution and $v_i$ is the average traffic speed of zone $i$. Clearly, \ref{30} indicates that the waiting time for Class 1 passengers is smaller than that of Class 2 passengers. This is intuitive as Class 1 passengers are more flexible and can enjoy a higher chance of getting matched to an idle vehicle regardless of the vehicle types.

We will employ logit models to capture the mode choices of passengers. Suppose  the percentage of Class 1 passengers out of the total population is $\theta\in [0,1]$, then the arrival rate of the two classes of passengers to the ride-hailing platform can be written as:
\begin{equation}
    \begin{cases}
    \lambda_{ij}^{AH} =\theta\lambda_{ij}^0  \frac{e^{-\epsilon c_{ij}^{AH}}}{e^{-\epsilon c_{ij}^{AH}}+e^{-\epsilon c_{ij}^0}};\\     
    \lambda_{ij}^{H} =(1-\theta)\lambda_{ij}^0  \frac{e^{-\psi c_{ij}^{H}}}{e^{-\psi c_{ij}^{H}}+e^{-\psi c_{ij}^0}},
    \end{cases}
    \label{2}\end{equation}
    
where $\lambda_{ij}^0$ denotes the arrival rate of potential travel demand from zone $i$ to zone $j$ regardless of the passenger classes and mode choices; $c_{ij}^0$ represents the cost of outside options (e.g., walking, transit, etc); $\epsilon$ and $\psi$ are sensitivity parameters of the logit models. Clearly, both $\lambda_{ij}^{AH}$ and $\lambda_{ij}^H$ are  decreasing functions with respect to the generalized travel cost of TNC trips, indicating that a higher TNC trip cost leads to the lower travel demand. 

\begin{remark}
The ride fare $r_i$ in equation (\ref{1}) only depends on the origin of the trip. This is consistent with the industry practice.  For instance, both Uber and Lyft  adopt a fixed per-time fare $p_t$ and per-distance fare $p_d$. The total trip fare is then calculated as the sum of a based fare, a time-based component and a distance-based component, multiplied by a surge multiplier $\gamma_i$ defined for each zone $i$ in order to motivate drivers to reposition to areas in short of supply. If we denote the average travel speed as $s$, then the total trip fare is proportional to $\gamma_i\left(p_t+s p_d\right)$, which depends on the trip origin only.
\end{remark}

\subsection{Incentives of the Mixed Fleet}
{\em Recruitment of AVs and human drivers}: The TNC platform recruits a mixture of AVs and for-hire human drivers to deliver mobility-on-demand services. Let $N_A$ and $N_H$ denote the total number of AVs and human-driver vehicles, respectively. Each AV is rented/purchased\footnote{The case of renting and purchasing does not make a fundamental difference for our model. Both will incur a cost for the platform, and their only difference is reflected in the value of this cost term.} from the manufacturer and incurs a cost, which includes both the procurement cost and the operational cost. On the other hand, each human driver is employed as an independent contractor, and the total number of human drivers depends on the hourly wage offered by the TNC platform. This can be captured by the following  supply function
\begin{equation}
\label{driversupply}
	{N_{H}} = {N_0}{F_d}(q),
\end{equation}
where $N_{0}$ denotes the total number of drivers who are looking for a job, $q$ denotes the average driver payment per unit time (e.g., an hour).\footnote{ The operational cost for human-driver vehicles, such as fuel, insurance, maintenance cost, are exogenous parameters that are not explicitly spelled out in (\ref{driversupply}). We comment that this does not incur any loss of generality as these exogenous parameters can be lumped into the driver supply function (\ref{driversupply}) as model parameters.} $F_d(q)$ represents the proportion of drivers who choose to work for the TNC platform. We assume that $F_d(q)$ is an increasing function, which follows the intuition that more drivers are willing to work for the platform if the net hourly wage $q$ is higher. Similar to the demand function, we do not assume $q$ to be the same for all drivers. In the rest of the paper, we assume that the human driver working choice (\ref{driversupply}) follows the logit model:
\begin{equation}
    N_H=N_0\cdot \frac{e^{\sigma q}}{e^{\sigma q}+e^{\sigma q_0}},
    \label{29}
\end{equation}
where $\sigma$ is the sensitivity parameter, and $q_0$ represents the earning of the outside options.

\begin{remark}
    In practice, the TNC platform cannot directly determine the hourly driver wage $q$. Instead, the platform can decide a commission rate $\delta$ that does not depend on the pickup and drop-off locations of each trip. Hence, the average driver wage $q$ is collectively determined by the commission rate $\delta$, passenger demand $\lambda_{ij}$, and the ride fare $r_i$ through
\begin{equation}
        q =(1-\delta) \dfrac{\displaystyle\sum_{i=1}^M \sum_{j=1}^M r_i \lambda_{ij}^{AH} t_{ij}}{N_H+N_A} + (1-\delta)\dfrac{\displaystyle\sum_{i=1}^M \sum_{j=1}^M r_i \lambda_{ij}^{H} t_{ij}}{N_H}.
    \end{equation}
Nonetheless, it is clear that $q$ and $\delta$ follow a one-to-one mapping. For this reason, we can replace $\delta$ with $q$ as a decision variable to simplify the solution and analysis without any loss of generality.
\end{remark}

{\em Characterization of driver status}: Each TNC vehicle will be in one of the three statuses: (1) delivering a passenger; (2) on the way to pick up a passenger; and (3) repositioning or cruising with empty seats to look for the next passenger. Let $N_{i,A}^I$ and $N_{i,H}^I$ denote number of idle AVs and human-driver vehicles in zone $i$, respectively, which corresponds to vehicles in the third status (including vehicles cruising in zone $i$, and passing by zone $i$ during repositioning). For the passengers without preference on AV or human drivers, we further assume that in each zone $i$, idle AVs and human drivers are randomly distributed, so that when the platform matches each passenger to the closest idle vehicle, there is a probability of $\dfrac{N_{i,A}^I}{N_{i,A}^I+N_{i,H}^I}$ to match to an AV, and a probability of  $\dfrac{N_{i,H}^I}{N_{i,A}^I+N_{i,H}^I}$ to match to a human-driver vehicle. In this case, the total vehicle hours supplied by AVs and human drivers can be decomposed into the sum of the following five terms\footnote{At equilibrium, we can treat total number of vehicles and the total supplied vehicles hours as equivalent.}:
\begin{align}
    \label{vehicle_conservation}
    \begin{cases}
        N_{A}=\displaystyle\sum_{i=1}^M\sum_{j=1}^M\dfrac{N_{i,A}^I}{N_{i,A}^I+N_{i,H}^I}\lambda_{ij}^{AH}t_{ij} \\
        +\displaystyle\sum_{i=1}^{M}\sum_{j=1}^{M}\dfrac{N_{i,A}^I}{N_{i,A}^I+N_{i,H}^I}\lambda_{ij}^{AH}w_i^{p,AH}
        +\sum_{i=1}^{M}N_{i,A}^I, \\
        N_{H}=\displaystyle\sum_{i=1}^{M}{\sum_{j=1}^M \left(\dfrac{N_{i,H}^I}{N_{i,A}^I+N_{i,H}^I}\lambda_{ij}^{AH}+\lambda_{ij}^{H}\right)t_{ij}} \\
        +\displaystyle\sum_{i=1}^{M}\sum_{j=1}^{M} \left({\dfrac{N_{i,H}^I}{N_{i,A}^I+N_{i,H}^I}\lambda_{ij}^{AH}w_i^{p,AH}}+\lambda_{ij}^{H}w_i^{p,H}\right)+\sum_{i=1}^{M}N_{i,H}^I.
    \end{cases}
\end{align}
where the first term in each equation of (\ref{vehicle_conservation}) denotes the hours of vehicles that is carrying a passenger, the second term represents the hours of vehicles that are on the way to pick up a passenger, and the third term corresponds to the hours that a vehicle is idle, including idle vehicles that are {\em both cruising in the same zone and repositioning to a different zone}.
\begin{remark}
We emphasize that the vehicles hours related to repositioning are already included in the third term of (\ref{vehicle_conservation}), thus we do not need to add it separately to the equation of vehicle conservation (\ref{vehicle_conservation}). This is a reasonable treatment because the repositioning vehicles can be matched to a nearby passenger at any time, even before it arrives at the repositioning destination. In other words, the repositioning vehicles are always open for matching as it passes by intermediate zones during the repositioning process, and therefore they should be counted in $N_{i,A}^I$ and $N_{i,H}^I$ when they pass by zone $i$. This treatment is different from some other works, such as \cite{li2019regulating}. However, such difference arises because the possibility of vehicles being intercepted during the repositioning process is not explicitly considered in \cite{li2019regulating}. In this case, the repositioning vehicles are not open for matching until they reach the destination, thus the repositioning hours can not be included in idle hours and it has to be counted as a separate term in (\ref{vehicle_conservation}). 
\end{remark}

The first two terms  in (\ref{vehicle_conservation}) are derived based on Little's law. If we denote $w_i^{d,A}$ and $w_i^{d,H}$ as the average vehicle waiting time between rides in zone $i$ for AVs and human-driven vehicles, respectively, then based on the same law, we can relate $w_i^{d,A}$ and $w_i^{d,H}$ to $N_{i,A}^I$ and $N_{i,H}^I$ through the following equations:
\begin{align}
\label{idle_term}
\begin{cases}
    N_{i,A}^I = w_i^{d,A}\cdot \displaystyle \sum_{j=1}^{M}\left( \dfrac{N_{i,A}^I}{N_{i,A}^I+N_{i,H}^I} \lambda_{ij}^{AH}\right),\\
    N_{i,H}^I = w_i^{d,H}\cdot \displaystyle\sum_{j=1}^{M}\left(\dfrac{N_{i,H}^I}{N_{i,A}^I+N_{i,H}^I} \lambda_{ij}^{AH} +  \lambda_{ij}^{H}\right).
\end{cases}
\end{align}

{\em Vehicle repositioning}: After each ride, the idle TNC vehicle can either cruise in the destination zone or strategically reposition itself to a different zone to seek the next passenger. For the TNC platform, the vehicle repositioning strategy is important as it can help reduce the spatial imbalance of supply and demand over the transportation network. In a mixed environment, the repositioning of AVs and human drivers are significantly different: AVs can be centrally dispatched by the platform to reposition between zones to maximize the platform profit, whereas human drivers selfishly and independently determine their own relocation strategies based on their perceived earning opportunities in distinct zones of the city\footnote{This is consistent with the current operational practices of Uber, Lyft, and Didi. }. Let $f_{ij,A}$ and $\Tilde{f}_{ij,H}$ denote the repositioning flow of AVs and human drivers, respectively, then $f_{ij,A}$ is a decision variable for the platform, while $\Tilde{f}_{ij,H}$ is an endogenous variable that depends on the driver's equilibrium choice, which will be modeled below. 

The human drivers' repositioning decisions depend on the expected income (per unit time) in each zone. For drivers in zone $i$, their average per-trip earning ${\bar{e}}_i$ correlates with the average trip time, the per-time trip fare, and the commission rate:
\begin{equation}
    {\bar{e}}_i = (1-\delta)  r_i {\bar{t}}_i = (1-\delta)  r_i  \cdot \frac{\displaystyle\sum_{j=1}^{M}t_{ij}{\left(\dfrac{ N_{i,H}^I}{N_{i,A}^I +N_{i,H}^I}\lambda_{ij}^{AH}+\lambda_{ij}^{H} \right)}}{\displaystyle\sum_{j=1}^{M}{\left(\dfrac{ N_{i,H}^I}{N_{i,A}^I +N_{i,H}^I}\lambda_{ij}^{AH}+ \lambda_{ij}^{H} \right)}},
    \label{6}
\end{equation}
where ${\bar{t}}_i$ represents the average duration of trips originating from zone $i$. To derive the per-time earning, we also need to characterize the average time spent by the driver to receive an earning of $\bar{e}_i$. After dropping off passengers in zone $i$, the human driver can either stay in zone $i$ or cruise to another zone $j$:
(a) if the driver chooses to stay in the same zone $i$, before completing the next ride, he/she will experience an average waiting time $w_j^{d,H}$ and an average trip time ${\bar{t}}_i$. The expected per-time earnings is then ${\bar{e}}_i \big/ (w_j^{d,H}+{\bar{t}}_i)$; (b) if the driver chooses to reposition to a distinct zone $j$, he/she will undergo an average waiting time $(w_j^{d,H}+t_{ij})$ and the average trip time ${\bar{t}}_j$, and thus the expected per-time earning is ${\bar{e}}_j \big/ (w_j^{d,H}+t_{ij}+{\bar{t}}_j)$. Given this, we can use a discrete choice model to characterize the probability that an idle human driver in zone $i$ chooses to reposition to zone $j$, which can be cast as:
\begin{align}
\begin{cases}
    \mathbb{P}_{ij}=\dfrac{e^{{\eta\bar{e}}_j/(w_j^{d,H}+t_{ij}+\bar{t}_j)}}{\sum_{k\neq i} e^{{\eta\bar{e}_k/(w_k^{d,H}+t_{ik}+\bar{t}_k)}}+e^{\eta\bar{e}_i/(w_j^{d,H}+\bar{t}_i)}}, \quad j\neq i, \\
    \mathbb{P}_{ii}=\dfrac{e^{{\eta\bar{e}}_j/(w_j^{d,H}+\bar{t}_j)}}{\sum_{k\neq i}e^{\eta\bar{e}_k/(w_k^{d,H}+t_{ik}+\bar{t}_k)}+e^{\eta\bar{e}_i/(w_j^{d,H}+\bar{t}_i)}}, \quad j = i.
\label{driver_reposition_probability}
\end{cases}
\end{align}
where $\mathbb{P}_{ij}$ is the probability that drivers in zone $i$ choose to reposition to zone $j$, and $\mathbb{P}_{ii}$ is the probability that an idle driver in zone $i$ chooses to stay in zone $i$. These repositioning probabilities are determined by the potential per-time earning in each zone. Given the repositioning probabilities, the intended repositioning flow of human drivers should satisfy:
\begin{equation}
    \label{repositioning_flow}
    f_{ij,H}=\mathbb{P}_{ij} \sum_{k=1}^M \left(\dfrac{ N_{k,H}^I}{N_{k,A}^I+N_{k,H}^I}\lambda_{ki}^{AH}+\lambda_{ki}^{H} \right),  \quad i,j\in\mathcal{V}.
\end{equation}
where the summation term on the right-hand side of (\ref{repositioning_flow}) denotes the total arrival rate of human drivers who drop off passengers in zone $i$ and need to decide whether to stay or to reposition to a different zone.

Note that $f_{ij,H}$ and $f_{ij,A}$ are the {\em intended} reposition flow before the vehicle starts repositioning. However, each vehicle can still be matched to a nearby passenger during the repositioning process. In other words, the empty vehicle can be intercepted by passengers in other zones before it arrives at the final repositioning destination. To model this, we denote $\mathcal{P}_{ij}\in\mathcal{V}$ as the set of intermediate zones along the shortest path from zone $i$ to zone $j$.\footnote{Note that $i\in\mathcal{P}_{ij}$, and $j\in\mathcal{P}_{ij}$.} Define $\Pi(k|i\rightarrow j)$ as the probability of the vehicle being intercepted when passing by $k$ while repositioning from zone $i$ to zone $j$. In this case,  the actual repositioning flow of AVs and human drivers can be written as:
\begin{align}
    \label{actualreposition}
    \begin{cases}
    \Tilde{f}_{ik,A} = \sum\limits_{j:k\in \mathcal{P}_{ij}}\Pi(k|i\rightarrow j) \cdot{f_{ij,A}},\quad i,j,k \in \mathcal{P}_{ij} \\
    \Tilde{f}_{ik,H} = \sum\limits_{j:k\in \mathcal{P}_{ij}}\Pi(k|i\rightarrow j) \cdot{f_{ij,H}},\quad i,j,k \in \mathcal{P}_{ij}
 
    \end{cases}
\end{align}
The interception probability $\Pi(k|i\rightarrow j)$ depends on the supply and demand in each zone along the path $\mathcal{P}_{ij}$, the order of passed zones in  $\mathcal{P}_{ij}$, and the consumed time to travel across these zones. We denote $d_i$ as the average time for passing by zone $i$ during repositioning, and denote $\varrho_i$ as the probability that the empty vehicle be matched with a nearby passenger before he/she leaves the zone $i$. To model $\varrho_i$, we consider an $M/M/1$ queue, where passengers are regarded as servers and drivers are regarded as jobs.\footnote{A similar queuing model has been used in \cite{banerjee2015pricing} to study optimal pricing strategies for ride-sharing platforms.} Based on this model, the waiting time $\tau$ of drivers for getting matched to a passenger in zone $i$ is subject to an exponential distribution \cite{harrison1993response}:
\begin{equation}
    \label{driverqueue}
    \mathbb{P}(\tau\leq T) = 1 - e^{-T/w_i^{d,H}}.
\end{equation}
Based on (\ref{driverqueue}), if an idle vehicle cruises in zone $i$ for a duration of $d_i$, then the probability of being intercepted by a passenger before it leaves zone $i$ is
\begin{equation}
    \label{matching proba}
    \varrho_i = \mathbb{P}(\tau\leq d_i) =  1 - e^{-d_i/w_i^{d,H}}.
\end{equation}
Now we can derive $\Pi(k|i\rightarrow j)$ based on $\varrho_i$. Specifically, if a vehicle is intercepted in zone $k$, it will not be intercepted in other zone in $\mathcal{P}_{ij}$ before entering zone $k$. For each of exposition, we order the set $\mathcal{P}_{ij}$ in an increasing sequence, i.e., the first element $p_1\in\mathcal{P}_{ij}$ is zone $i$; the second  element $p_2\in\mathcal{P}_{ij}$ is the closest zone to zone $i$ along the path $\mathcal{P}_{ij}$; the third element $p_3\in\mathcal{P}_{ij}$ is the closest zone to $p_2$ along the path $\mathcal{P}_{ij}$, etc. Let $k$ be the $s$th element, i.e., $p_s=k$, then for a specific $k\neq j$, the interception probability  $\Pi(k|i\rightarrow j)$ is modelled as:
\begin{equation}
    \label{interceptionprobability}
    \Pi(k|i\rightarrow j) = (1-\varrho_{p_1})\cdot(1-\varrho_{p_2})\cdots(1-\varrho_{p_{s-1}})\cdot\varrho_{p_s}.
\end{equation}
Equation (\ref{interceptionprobability}) indicates that if the repositioning vehicle is intercepted at zone $p_s = k$, it must have travelled across the zones $p_1, p_2, \cdots, p_{s-1}$ without being intercepted before entering zone $k$.


{\em Flow balancing:} In the long term, for each zone $i$, the inflow and outflow of AVs (or human drivers) should be equal. This yields the following balancing constraints:
\begin{subnumcases}
{\label{split_flow_rebalance}}
    \sum_{j=1}^{M}\left({\dfrac{N_{j,A}^I}{N_{j,A}^I+N_{j,H}^I}}\lambda_{ji}^{AH}+\tilde{f}_{ji,A}\right)=\sum_{j=1}^{M}\left(\dfrac{N_{i,A}^I}{N_{i,A}^I+N_{i,H}^I}\lambda_{ij}^{AH}+\tilde{f}_{ij,A}\right),\forall i\in\mathcal{V},
    \label{split_flow_rebalance_A}\\
    \sum_{j=1}^{M}\left(\dfrac{N_{j,H}^I}{N_{j,A}^I+N_{j,H}^I}\lambda_{ji}^{AH}+\lambda_{ji}^{H}+\Tilde{f}_{ji,H}\right)=\sum_{j=1}^{M}\left(\dfrac{N_{i,H}^I}{N_{i,A}^I+N_{i,H}^I}\lambda_{ij}^{AH}+\lambda_{ij}^{H}+\Tilde{f}_{ij,H}\right),\forall i\in\mathcal{V}.
\end{subnumcases}
where the left-hand side denotes the total number of AVs (or human drivers) that flow into zone $i$ per unit time (including both occupied vehicles and idle repositioning vehicles), and the right-hand side represents the total number of AVs (or human drivers) that exit from zone $i$ and flow into other zones. They must be balanced at the equilibrium.

\emph{Traffic congestion:} The ride-hailing vehicles will inevitably affect the traffic congestion of the city. In general, these vehicles contribute to the overall traffic speed that has a great impact on the trip time from zone $i$ to $j$. In this case, the trip time $t_{ij}$ should be an endogenous variable that relies on the spatial distribution of ride-hailing vehicles. This work propose a traffic congestion model to model the endogenous relationship. In general, the ride-hailing vehicles are distributed following an asymmetric pattern over the city. In the urban area of the city, the streets are usually with heavy congestion, launching a lot of ride-sourcing vehicles will slow down the traffic speed. On contrary, in remote areas, the traffic congestion rarely happens so that the ride-hailing vehicles have a negligible effects on the overall traffic speed. To model the endogenous relationship, we employ notation $\mathcal{C} \in \mathcal{V}$ and  $\mathcal{R} \in \mathcal{V}\setminus\mathcal{C}$ to represent all the ``congested urban zones''  and the ``uncongested remote zones'' in a city, respectively. Moreover, we use notation $v_c$ and $v_r$ to indicate the traffic speeds in urban and remote areas, respectively, so that $v_i=v_c$ if $i\in \mathcal{C}$ and $v_i=v_r$ if $i\in \mathcal{R}$. Further, let $d_{ij}$ represent the shortest distance from zone $i$ to $j$. A trip may travel across the two areas  $\mathcal{C}$ and $\mathcal{R}$, then the shortest path $d_{ij}$ is comprised of $d_{ij}^\mathcal{C}$ and $d_{ij}^\mathcal{R}$, which res presents the distance of the shortest path that lies within $\mathcal{C}$ and $\mathcal{R}$, respectively. The travel time $t_{ij}$ can be then derived as:
\begin{equation}
    \label{traveltime}
    t_{ij} = d_{ij}^\mathcal{C}\cdot \frac{1}{v_c} + d_{ij}^\mathcal{R}\cdot \frac{1}{v_r}.
\end{equation}
where the first term denotes the trip time in congestion zones $\mathcal{C}$ and the second term represents the trip time in uncongested zones $\mathcal{R}$. We further assume the traffic speed in remote zones is not affected by the ride-hailing vehicle fleet. In this case, the parameters $d_{ij}^\mathcal{C}, d_{ij}^\mathcal{R}$ and $v_r$ are all exogenous parameters in equation (\ref{traveltime}). Conversely, the traffic in congested zones is highly correlated with the ride-hailing vehicle fleets. Let us use $N_\mathcal{C}$ to represent the total number of vehicles in congested zones $\mathcal{C}$. Intuitively, the overall traffic speed  $v_c$ in $\mathcal{C}$ is a decreasing function $N_\mathcal{C}$. Therefore, the trip time $t_{ij}$ is also a function of $N_\mathcal{C}$.

Finally, the total number of ride-hailing vehicle fleets in $\mathcal{C}$ should encompass all the vehicles in congested zones $\mathcal{C}$, including occupied vehicles, cruising vehicles, and empty vehicles on way to pick up passengers. Inspired by equation (\ref{vehicle_conservation}), $N_\mathcal{C}$ can be calculated as:
\begin{equation}
    \label{congestvehicle}
    N_\mathcal{C} = \displaystyle\sum_{i=1}^{M}{\sum_{j=1}^M\left({\lambda_{ij}^{AH}+\lambda_{ij}^{H}}\right)}d_{ij}^\mathcal{C}\cdot \frac{1}{v_c}+\sum_{i\in \mathcal{C}}\sum_{j=1}^{M}\left({\lambda_{ij}^{AH}w_i^{p,AH}+\lambda_{ij}^{H}w_i^{p,H}}\right)+\sum_{i\in \mathcal{C}}{\left(N_{i,A}^I+N_{i,H}^I\right)}.
\end{equation}
where the first term is trip time in congestion zones $\mathcal{C}$; the second term denotes all vehicles on the way to pick up the passenger in $\mathcal{C}$; the third term denotes all idle vehicles in $\mathcal{C}$. Each term in equation (\ref{congestvehicle}) reflects a proportion of the corresponding terms in equation (\ref{vehicle_conservation}), which only corresponds to the congested zones $\mathcal{C}$ with mixed fleets.
\subsection{Incentives of TNC Platform}
The TNC platform hires the mixed fleet and makes operational decisions to maximize its profit. In this case, the decisions of the platform include the ride fare $r_i$, human driver wage $q$, AV fleet size $N_A$, the spatial distribution of idle autonomous vehicles $N_{i,A}^I$, and the repositioning flow of autonomous vehicles $f_{ij,A}$. The platform receives revenue from the passenger and bears the cost incurred by AVs and human drivers. Let ${\bf r}=(r_1, \ldots, r_M)$, ${\bf N}_{A}^I=(N_{1,A}^I, \ldots, N_{M,A}^I)$, and denote the AV repositioning flow as a vector ${\bf f}_A$, then the profit maximization problem of the TNC platform can be written as:
\begin{equation}
\hspace{-6cm}
\label{optimalpricing_trip}
    \mathop {\max }\limits_{ q, {\bf r}, {\bf N}_{A}^I, {\bf f}_A }\sum_{i=1}^{M}\sum_{j=1}^{M}{r_i\left(\lambda_{ij}^{AH}+\lambda_{ij}^{H}\right)t_{ij}}-C\cdot N_{A}-N_0F_d\left(q\right)q
\end{equation}
\begin{subnumcases} {\label{profit_constraints}}
    \lambda_{ij}^{AH} =\theta\lambda_{ij}^0  \frac{e^{-\epsilon c_{ij}^{AH}}}{e^{-\epsilon c_{ij}^{AH}}+e^{-\epsilon c_{ij}^0}}, \label{10a}\\
    \lambda_{ij}^{H} =(1-\theta)\lambda_{ij}^0  \frac{e^{-\psi c_{ij}^{H}}}{e^{-\psi c_{ij}^{H}}+e^{-\psi c_{ij}^0}}, \label{b}\\
    w_i^{p,AH}(N_i^I) \leq w_i^{p,max}, \quad w_i^{p,A}(N_i^{I,A}) \leq w_i^{p,max}, \label{10c}\\
    N_{i,A}^I = w_i^{d,A}\cdot \displaystyle \sum_{j=1}^{M}\left( \dfrac{N_{i,A}^I}{N_{i,A}^I+N_{i,H}^I} \lambda_{ij}^{AH}\right),\label{10dd} \\
    N_{i,H}^I = w_i^{d,H}\cdot \displaystyle\sum_{j=1}^{M}\left(\dfrac{N_{i,H}^I}{N_{i,A}^I+N_{i,H}^I} \lambda_{ij}^{AH} +  \lambda_{ij}^{H}\right), \label{10d}\\
    \sum_{j=1}^{M}\left({\dfrac{N_{j,A}^I}{N_{j,A}^I+N_{j,H}^I}}\lambda_{ji}^{AH}+\tilde{f}_{ji,A}\right)=\sum_{j=1}^{M}\left(\dfrac{N_{i,A}^I}{N_{i,A}^I+N_{i,H}^I}\lambda_{ij}^{AH}+\tilde{f}_{ij,A}\right),\forall i\in\mathcal{V},\label{10f}\\
    \sum_{j=1}^{M}\left(\dfrac{N_{j,H}^I}{N_{j,A}^I+N_{j,H}^I}\lambda_{ji}^{AH}+\lambda_{ji}^{H}+\Tilde{f}_{ji,H}\right)=\sum_{j=1}^{M}\left(\dfrac{N_{i,H}^I}{N_{i,A}^I+N_{i,H}^I}\lambda_{ij}^{AH}+\lambda_{ij}^{H}+\Tilde{f}_{ij,H}\right),\forall i\in\mathcal{V}, \label{10g}\\
    \label{10h}{N_0}{F_d}(q)=\displaystyle\sum_{i=1}^{M}{\sum_{j=1}^M \left(\dfrac{N_{i,H}^I}{N_{i,A}^I+N_{i,H}^I}\lambda_{ij}^{AH}+\lambda_{ij}^{H}\right)t_{ij}} \\ \notag
    \quad \quad \quad \quad \quad +\displaystyle\sum_{i=1}^{M}\sum_{j=1}^{M} \left({\dfrac{N_{i,H}^I}{N_{i,A}^I+N_{i,H}^I}\lambda_{ij}^{AH}w_i^{p,AH}}+\lambda_{ij}^{H}w_i^{p,H}\right)+\sum_{i=1}^{M}N_{i,H}^I,\\\label{10i}
    N_{A}=\displaystyle\sum_{i=1}^M\sum_{j=1}^M\dfrac{N_{i,A}^I}{N_{i,A}^I+N_{i,H}^I}\lambda_{ij}^{AH}t_{ij} 
    +\displaystyle\sum_{i=1}^{M}\sum_{j=1}^{M}\dfrac{N_{i,A}^I}{N_{i,A}^I+N_{i,H}^I}\lambda_{ij}^{AH}w_i^{p,AH}
    +\sum_{i=1}^{M}N_{i,A}^I,\\
     N_\mathcal{C} = \displaystyle\sum_{i=1}^{M}{\sum_{j=1}^M\left({\lambda_{ij}^{AH}+\lambda_{ij}^{H}}\right)}d_{ij}^\mathcal{C}\cdot \frac{1}{v_c}+\sum_{i\in \mathcal{C}}\sum_{j=1}^{M}\left({\lambda_{ij}^{AH}w_i^{p,AH}+\lambda_{ij}^{H}w_i^{p,H}}\right)+\sum_{i\in \mathcal{C}}{\left(N_{i,A}^I+N_{i,H}^I\right)},\label{10j}
\end{subnumcases}
where $N_A$ is the total number of AV deployed by platform, $C$ is capital and operational cost for each AV per unit time\footnote{We implicitly assume that the operational cost of AVs does not depend on the status of the vehicle. This is a reasonable assumption if the travel speed of AVs does not change significantly across different vehicles statuses.}, and $\mathbb{P}_{ij}$ is determined by the logit model (\ref{driver_reposition_probability}). The objective function (\ref{optimalpricing_trip}) denotes the platform profit. The first term \textcolor{blue}{$\sum_{i=1}^{M}\sum_{j=1}^{M}{r_i\left(\lambda_{ij}^{AH}+\lambda_{ij}^{H}\right)t_{ij}}$} in (\ref{optimalpricing_trip}) defines the total revenue, the second term $C \cdot N_{A}$ represents the AV deployment cost, and the last term $N_0F_d\left(q\right)q$ denotes the total payment to human drivers. Constraints (\ref{10a}) and (\ref{b}) showcase the demand of passengers. Constraint (\ref{10c}) imposes an upper bound of passenger waiting time\footnote{In specific regions of the urban area, the demand for transportation is remarkably low, thereby making it economically rational for service providers to refrain from offering any services to these zones. To prevent the emergence of trivial solutions resulting from this scenario, we enforce an upper bound on the passenger waiting time.}. Constraints (\ref{10dd}) and (\ref{10d}) follow equation (\ref{idle_term}).  Equations (\ref{10f}) and (\ref{10g}) are the flow balance constraints in each zone for two types of vehicles following (\ref{split_flow_rebalance}), where $\tilde{f}_{ij,A}$ and $\tilde{f}_{ij,H}$ satisfy (\ref{driver_reposition_probability})-(\ref{interceptionprobability}). 
Constraint (\ref{10h}) combines the driver supply model (\ref{driversupply}) and the human-driven vehicle hour conservation constraint (\ref{vehicle_conservation}). Constraint (\ref{10i}) and (\ref{10j}) denote the conservation of vehicles. We observe that the above optimization problem is non-concave. Formally, we have the following lemma:
\begin{lemma} The profit-maximization problem of the TNC platform (\ref{optimalpricing_trip})-(\ref{profit_constraints}) is non-concave.
\end{lemma}
\begin{proof}  For the optimization problem (\ref{optimalpricing_trip})-(\ref{profit_constraints}) to be concave, the objective function (\ref{optimalpricing_trip}) should be jointly concave with respect to the decision variables $(q, {\bf r}, {\bf N}_{A}^I, {\bf f}_A)$, and concurrently,  the constraints (\ref{profit_constraints}) should also be linear with respect to $(q, {\bf r}, {\bf N}_{A}^I, {\bf f}_A)$ \cite{boyd2004convex}. Regardless of whether the objective function is concave or not, it is clear that all the constraints in (\ref{profit_constraints}) are  nonlinear (especially in the repositioning relations of human-driven vehicles). Therefore, we can conclude that the profit-maximization problem of the TNC platform (\ref{optimalpricing_trip})-(\ref{profit_constraints}) is non-concave.
\end{proof}

\begin{remark} In this paper, we assume that AVs
can be deployed in the entire study area. This is not an unreasonable assumption for many
practical cases. For instance, currently Waymo and Cruise are operating in the entire San
Francisco city. In this case, the ``remote areas” can be understood as suburban areas and/or lower-density urban areas (e.g., the west and south of San Francisco city), where AVs can also be deployed. However, when rural areas are considered, the nature of the problem can be different because  additional infrastructure (e.g., high-
definition map, road-side V2X infrastructure, etc) may be needed to support the deployment
of AVs. In this case, a more complicated transportation network should be considered, where part
of the study area is AV-enabling and the rest is not. This is an interesting research direction,
but is beyond the scope of our paper. 
\end{remark}

\begin{remark} In the ideal case, the modelling of traffic congestion would occur at the zonal level - each specific area within the urban landscape would possess its own traffic velocity, contingent upon the aggregate quantity of ride-sourcing vehicles present within the given zone. However, our model does not accommodate for this degree of granularity in capturing congestion data due to the substantial amount of requisite data (e.g., background traffic, roadway capacity), which remains beyond our current access. As such, the integration of a more comprehensive congestion model is reserved for future research endeavours. \end{remark}

\section{The Solution Algorithm}
\label{The solution algorithm}
The optimization problem (\ref{optimalpricing_trip}) is non-concave due to the complex non-linearity in the demand function, the supply function, the passenger waiting time model, and the human-driver repositioning model. Therefore, it is challenging to numerically derive the globally optimal solution to  (\ref{optimalpricing_trip}). In this section, we will present an algorithm that provides an approximate solution to (\ref{optimalpricing_trip}), and at the same time we will develop a decomposition-based method to rigorously quantify the distance between our proposed solution and the unknown globally optimal solution.
\subsection{Constraint relaxation and upper bound}
We note that the complex coupling between distinct zones in the transportation network primarily arises from the repositioning flow, which are captured by constraints (\ref{10f}) and (\ref{10g}). To obtain an approximate solution to (\ref{optimalpricing_trip}), we propose to relax (\ref{10f}) and (\ref{10g}), and then solve the relaxed problem exactly, which offers 
an upper bound on the quality of solution compared to that of the original problem. Note that without (\ref{10f}) and (\ref{10g}), constraints (\ref{10dd}) and (\ref{10d})  can be also removed without loss of generality because they are only relevant if (\ref{10f}) is considered. In this case, the relaxed problem can be formulated as:
\begin{equation}
\hspace{-6cm}
\label{optimalpricing_trip_relax}
    \mathop {\max }\limits_{ q, {\bf r}, {\bf N}_{A}^I, {\bf N}_{H}^I }\sum_{i=1}^{M}\sum_{j=1}^{M}{r_i\left(\lambda_{ij}^{AH}+\lambda_{ij}^{H}\right)t_{ij}}-C\cdot N_{A}-N_0F_d\left(q\right)q.
\end{equation}
\begin{subnumcases}
   {\label{profit_constraints_relax}}
   \lambda_{ij}^{AH} =\theta\lambda_{ij}^0  \frac{e^{-\epsilon c_{ij}^{AH}}}{e^{-\epsilon c_{ij}^{AH}}+e^{-\epsilon c_{ij}^0}}, \label{14a}\\
    \lambda_{ij}^{H} = (1-\theta)\lambda_{ij}^0  \frac{e^{-\psi c_{ij}^{H}}}{e^{-\psi c_{ij}^{H}}+e^{-\psi c_{ij}^0}}, \label{14a2}\\
    w_i^{p,AH}(N_i^I) \leq w_i^{p,max}, \quad w_i^{p,A}(N_i^{I,A}) \leq w_i^{p,max},\label{14aa}\\
    \label{14d} {N_0}{F_d}(q)=\displaystyle\sum_{i=1}^{M}{\sum_{j=1}^M \left(\dfrac{N_{i,H}^I}{N_{i,A}^I+N_{i,H}^I}\lambda_{ij}^{AH}+\lambda_{ij}^{H}\right)t_{ij}} \\ \notag
    \quad \quad \quad \quad \quad +\displaystyle\sum_{i=1}^{M}\sum_{j=1}^{M} \left({\dfrac{N_{i,H}^I}{N_{i,A}^I+N_{i,H}^I}\lambda_{ij}^{AH}w_i^{p,AH}}+\lambda_{ij}^{H}w_i^{p,H}\right)+\sum_{i=1}^{M}N_{i,H}^I, \\
    N_{A}=\displaystyle\sum_{i=1}^M\sum_{j=1}^M\dfrac{N_{i,A}^I}{N_{i,A}^I+N_{i,H}^I}\lambda_{ij}^{AH}t_{ij} 
    +\displaystyle\sum_{i=1}^{M}\sum_{j=1}^{M}\dfrac{N_{i,A}^I}{N_{i,A}^I+N_{i,H}^I}\lambda_{ij}^{AH}w_i^{p,AH}
    +\sum_{i=1}^{M}N_{i,A}^I,
    \label{14e}\\
    N_\mathcal{C} =\displaystyle\sum_{i=1}^{M}{\sum_{j=1}^M\left({\lambda_{ij}^{AH}+\lambda_{ij}^{H}}\right)}d_{ij}^\mathcal{C}\cdot \frac{1}{v_c}+\sum_{i\in \mathcal{C}}\sum_{j=1}^{M}\left({\lambda_{ij}^{AH}w_i^{p,AH}+\lambda_{ij}^{H}w_i^{p,H}}\right)+\sum_{i\in \mathcal{C}}{\left(N_{i,A}^I+N_{i,H}^I\right)}.
    \label{14f}
\end{subnumcases}
where ${\bf N}_{H}^I=(N_{1,H}^I,\ldots,N_{H,H}^I)$ is the vector that denotes the number of idle human drivers in each zone.
\subsection{The solution algorithm}
Prior to presenting the algorithm's details, it is necessary to reconsider constraint (\ref{14f}), which involves $N_\mathcal{C}$ as a fixed point. We note that when  $N_\mathcal{C}$ is given, the traffic speed $v_c$ and trip time $t_{ij}$ are uniquely determined. Consequently, the problem in question can be resolved by iterating over $N_\mathcal{C}$ and subsequently conducting a dual decomposition for each value of $N_\mathcal{C}$. Specifically, we begin by carrying out a grid search on the scalar decision variable $N_\mathcal{C}$, which allows us to generate the corresponding $v_c$ and $t_{ij}$ values for each $N_\mathcal{C}$. By plugging (\ref{14e}) into (\ref{optimalpricing_trip_relax}), the objective function becomes separable over each zone. The only complicating constraints that couple the decision variables of different zones are (\ref{14d}) and (\ref{14f}). To address this constraint, we form the separable Lagrangian:
\begin{multline} \label{separable lagrangian}
    L(q, {\bf r}, {\bf N}_{A}^I, {\bf N}_{H}^I) = \sum_{i=1}^{M}\sum_{j=1}^{M}{r_i\left(\lambda_{ij}^{AH}+\lambda_{ij}^{H}\right)t_{ij}}
    - C\left( \displaystyle\sum_{i=1}^M\sum_{j=1}^M \dfrac{N_{i,A}^I}{N_{i,A}^I+N_{i,H}^I}\lambda_{ij}^{AH}t_{ij}\right) \\
    - C\left(\displaystyle\sum_{i=1}^{M}\sum_{j=1}^{M} \dfrac{N_{i,A}^I}{N_{i,A}^I+N_{i,H}^I}\lambda_{ij}^{AH}w_i^{p,AH}
    +\sum_{i=1}^{M}N_{i,A}^I \right) \\
    +\mu \left({N_0}{F_d}(q)-\displaystyle\sum_{i=1}^{M}{\sum_{j=1}^M \dfrac{N_{i,H}^I}{N_{i,A}^I+N_{i,H}^I}\lambda_{ij}^{AH}t_{ij}}-\displaystyle\sum_{i=1}^M\sum_{j=1}^M \lambda_{ij}^{H}t_{ij}\right)\\
    +\mu \left(-\displaystyle\sum_{i=1}^{M}\sum_{j=1}^{M} {\dfrac{N_{i,H}^I}{N_{i,A}^I+N_{i,H}^I}\lambda_{ij}^{AH}w_i^{p,AH}}-\sum_{i=1}^{M}\sum_{j=1}^{M}\lambda_{ij}^{H}w_i^{p,H}-\sum_{i=1}^{M}N_{i,H}^I\right)\\
    +\kappa\left( \displaystyle\sum_{i=1}^{M}{\sum_{j=1}^M\left({\lambda_{ij}^{AH}+\lambda_{ij}^{H}}\right)}d_{ij}^\mathcal{C}\cdot \frac{1}{v_c}+\sum_{i\in \mathcal{C}}\sum_{j=1}^{M}\left({\lambda_{ij}^{AH}w_i^{p,AH}+\lambda_{ij}^{H}w_i^{p,H}}\right)  +\sum_{i\in \mathcal{C}}{\left(N_{i,A}^I+N_{i,H}^I\right)}\right)\\
    -N_0F_d\left(q\right)q,
\end{multline}
where $\mu,\kappa$ are the Lagrange multipliers associated with (\ref{14d}) and (\ref{14f}), respectively. Note that $\lambda_{ij}^{AH}$ and $\lambda_{ij}^H$ are  functions of $r_i$, $N_{i,A}^I$ and $N_{i,H}^I$. Clearly, (\ref{separable lagrangian}) can be decomposed into $M+1$ terms, i.e.,
\[L(q, {\bf r}, {\bf N}_{A}^I, {\bf N}_{H}^I) = \sum_{i=1}^M L_i( r_i, {N}_{i,A}^I, {N}_{i,H}^I)+ L_q(q),\]
in which the term corresponding to the decisions in zone $i$ can be written down in two cases:\\
(1): if $i \in \mathcal{C}$:
\begin{align}
\label{decomposed_i}
    &L_i( r_i, {N}_{i,A}^I, {N}_{i,H}^I)=L_i^{\mathcal{C}}( r_i, {N}_{i,A}^I, {N}_{i,H}^I)= 
     \sum_{j=1}^{M}{r_i\left(\lambda_{ij}^{AH}+\lambda_{ij}^{H}\right)t_{ij}}\nonumber\\
    &-C\cdot \left( {\sum_{j=1}^{M}{\dfrac{N_{i,A}^I}{N_{i,A}^I+N_{i,H}^I}\lambda_{ij}^{AH}t_{ij}}}+\sum_{j=1}^{M}\dfrac{N_{i,A}^I}{N_{i,A}^I+N_{i,H}^I}\lambda_{ij}^{AH}w_i^{p,AH}+N_{i,A}^I \right)\nonumber\\
    & -\mu \cdot\left( {\sum_{j=1}^{M} {\dfrac{N_{i,H}^I}{N_{i,A}^I+N_{i,H}^I}\lambda_{ij}^{AH}t_{ij}}}+\sum_{j=1}^{M}{\lambda_{ij}^{H}t_{ij}}+\sum_{j=1}^{M}{\dfrac{N_{i,H}^I}{N_{i,A}^I+N_{i,H}^I}\lambda_{ij}^{AH}w_i^{p,AH}}+\sum_{j=1}^M \lambda_{ij}^{H}w_i^{p,H}+N_{i,H}^I\right)\\ \notag
    & - \kappa\cdot\left( \displaystyle{\sum_{j=1}^M\left({\lambda_{ij}^{AH}+\lambda_{ij}^{H}}\right)}d_{ij}^\mathcal{C}\cdot \frac{1}{v_c}+\sum_{j=1}^{M}{\left({\lambda_{ij}^{AH}w_i^{p,AH}+\lambda_{ij}^{H}w_i^{p,H}}\right)}+{\left(N_{i,A}^I+N_{i,H}^I\right)}\right),
\end{align}
(2): if $i \in \mathcal{R}$:
\begin{align}
\label{decomposed_i2}
   & L_i( r_i, {N}_{i,A}^I, {N}_{i,H}^I)=L_i^{\mathcal{R}}( r_i, {N}_{i,A}^I, {N}_{i,H}^I)= \sum_{j=1}^{M}{r_i\left(\lambda_{ij}^{AH}+\lambda_{ij}^{H}\right)t_{ij}} \nonumber\\
    &-C\cdot \left( {\sum_{j=1}^{M}{\dfrac{N_{i,A}^I}{N_{i,A}^I+N_{i,H}^I}\lambda_{ij}^{AH}t_{ij}}}+\sum_{j=1}^{M}\dfrac{N_{i,A}^I}{N_{i,A}^I+N_{i,H}^I}\lambda_{ij}^{AH}w_i^{p,AH}+N_{i,A}^I \right)\nonumber\\
    & -\mu \cdot\left( {\sum_{j=1}^{M} {\dfrac{N_{i,H}^I}{N_{i,A}^I+N_{i,H}^I}\lambda_{ij}^{AH}t_{ij}}}+\sum_{j=1}^{M}{\lambda_{ij}^{H}t_{ij}}+\sum_{j=1}^{M}{\dfrac{N_{i,H}^I}{N_{i,A}^I+N_{i,H}^I}\lambda_{ij}^{AH}w_i^{p,AH}}+\sum_{j=1}^M \lambda_{ij}^{H}w_i^{p,H}+N_{i,H}^I\right)\\ \notag
    & - \kappa\cdot\displaystyle{\sum_{j=1}^M\left({\lambda_{ij}^{AH}+\lambda_{ij}^{H}}\right)}d_{ij}^\mathcal{C}\cdot \frac{1}{v_c},
\end{align}
and the term corresponding to $q$ is
\begin{equation}
    \label{decomposed_q}
    L_q(q) = \mu N_0F_d(q) - N_0F_d(q)q.
\end{equation}
Therefore, given the Lagrange multipliers $\mu$ and $\kappa$, the relaxed problem (\ref{optimalpricing_trip_relax}) can be decomposed into $M+1$ small-scale optimization problem, including
\begin{equation}
    \hspace{-4cm}
\label{decomposed_opt_i}
    \mathop {\max }\limits_{ r_i\geq 0, {N}_{i,A}^I\geq 0, {N}_{i,H}^I\geq 0 }L_i( r_i, {N}_{i,A}^I, {N}_{i,H}^I)
\end{equation}
\begin{subnumcases}
   {\label{profit_constraints_relax4}}
    \lambda_{ij}^{AH} =\theta\lambda_{ij}^0  \frac{e^{-\epsilon c_{ij}^{AH}}}{e^{-\epsilon c_{ij}^{AH}}+e^{-\epsilon c_{ij}^0}},
    \label{16a}\\
    \lambda_{ij}^{H} = (1-\theta)\lambda_{ij}^0  \frac{e^{-\psi c_{ij}^{H}}}{e^{-\psi c_{ij}^{H}}+e^{-\psi c_{ij}^0}}, \label{16c}\\
    w_i^{p,AH}(N_i^I) \leq w_i^{p,max}, \quad w_i^{p,A}(N_i^{I,A}) \leq w_i^{p,max},\label{16b}
\end{subnumcases}
and
\begin{equation}
    \hspace{-2.5cm}
\label{decomposed_opt_q}
    \mathop {\max }\limits_{ q\geq 0}L_q(q).
\end{equation}
The decomposition procedure described above enables us to develop a dual-decomposition algorithm that efficiently computes the globally optimal solution to (\ref{optimalpricing_trip_relax}). The pseudo-code of the algorithm is summarized in Algorithm \ref{algorithm1}.

\begin{algorithm}[htbp]
\caption{The solution algorithm for the optimal spatial pricing problem (\ref{optimalpricing_trip})}\label{algorithm1}
\begin{algorithmic}[1]
\REQUIRE Initial guess of the primal variables $(\bar{q}, {\bf \bar{r}}, {\bf \bar{N}}_{A}^I, {\bf \bar{N}}_{H}^I)$, and  the dual variables $\mu, \kappa$.
\STATE Setup stopping criterion: constraints (\ref{14d})-(\ref{14e})  satisfied or maximum iteration reached. 
\FOR {each iteration}
\STATE {Enumerate over $N_\mathcal{C}$ by selecting and fixing the value of $N_\mathcal{C}$ in each iteration}
\STATE {Calculate $v_c$ and $t_{ij}$ based on the selected $N_\mathcal{C}$ in each iteration,}
\WHILE {stopping criterion not satisfied}
\STATE {{For each $i\in \mathcal{C}$, solve the sub-problem for zone $i$:}
\begin{equation}
    \hspace{-4cm}
\label{primal_dual_i}
    \mathop {\max }\limits_{ r_i\geq 0, {N}_{i,A}^I\geq 0, {N}_{i,H}^I\geq 0 } L_i^{\mathcal{C}}( r_i, {N}_{i,A}^I, {N}_{i,H}^I)
\end{equation}
\begin{subnumcases}
   {\label{primal_dual_constinat1}}
    \lambda_{ij}^{AH} =\theta\lambda_{ij}^0  \frac{e^{-\epsilon c_{ij}^{AH}}}{e^{-\epsilon c_{ij}^{AH}}+e^{-\epsilon c_{ij}^0}},\\
    \lambda_{ij}^{H} = (1-\theta)\lambda_{ij}^0  \frac{e^{-\psi c_{ij}^{H}}}{e^{-\psi c_{ij}^{H}}+e^{-\psi c_{ij}^0}},\\
   w_i^{p,AH}(N_i^I) \leq w_i^{p,max}, \quad w_i^{p,A}(N_i^{I,A}) \leq w_i^{p,max},
\end{subnumcases}
}
\STATE {For each $i\in \mathcal{R}$, solve the sub-problem for zone $i$:}
\begin{equation}
    \hspace{-1.5cm}
    \label{primal_dual_ii}
    \mathop {\max }\limits_{ r_i\geq 0, {N}_{i,A}^I\geq 0, {N}_{i,H}^I\geq 0 } L_i^{\mathcal{R}}( r_i, {N}_{i,A}^I, {N}_{i,H}^I)\end{equation}
    \centerline{Subject to (\ref{primal_dual_constinat1})}
\STATE Solve the decomposed Lagrangian over driver supply,
\begin{equation}
    \hspace{-2.5cm}
\label{primal_dual_q}
    \mathop {\max }\limits_{ q\geq 0}L_q(q)
\end{equation}
\STATE Choose an appropriate step sizes $\tau_1, \tau_2$ and update the dual variables $\mu,\kappa$,
{
\begin{subnumcases}
{\label{primaldual3}}
\mu = \mu - \tau_1 \cdot \left(    {N_0}{F_d}(q)-\displaystyle\sum_{i=1}^{M}{\sum_{j=1}^M \dfrac{N_{i,H}^I}{N_{i,A}^I+N_{i,H}^I}\lambda_{ij}^{AH}t_{ij}}-\displaystyle\sum_{i=1}^M\sum_{j=1}^M \lambda_{ij}^{H}t_{ij}\right)\\\notag
\quad \quad +\tau_1\left(-\displaystyle\sum_{i=1}^{M}\sum_{j=1}^{M} {\dfrac{N_{i,H}^I}{N_{i,A}^I+N_{i,H}^I}\lambda_{ij}^{AH}w_i^{p,AH}}-\sum_{i=1}^{M}\sum_{j=1}^{M}\lambda_{ij}^{H}w_i^{p,H}-\sum_{i=1}^{M}N_{i,H}^I\right), \\
\kappa = \kappa - \tau_2 \cdot \left( N_\mathcal{C} - \displaystyle\sum_{i=1}^{M}{\sum_{j=1}^M\left({\lambda_{ij}^{AH}+\lambda_{ij}^{H}}\right)}d_{ij}^\mathcal{C}\cdot \frac{1}{v_c}-\sum_{i\in \mathcal{C}}\sum_{j=1}^{M}{\left({\lambda_{ij}^{AH}w_i^{p,AH}+\lambda_{ij}^{H}w_i^{p,H}}\right)}\right)\\ \notag
\quad \quad -\tau_2 \cdot \left(- \sum_{i\in \mathcal{C}}{\left(N_{i,A}^I+N_{i,H}^I\right)}\right),
\end{subnumcases}}
\ENDWHILE 
\ENDFOR
\STATE Obtain the solution $(\bar{q}, {\bf \bar{r}}, {\bf \bar{N}}_{A}^I, {\bf \bar{N}}_{H}^I)$ and the corresponding platform profit $\bar{R}$. 
\STATE Use $(\bar{q}, {\bf \bar{r}}, {\bf \bar{N}}_{A}^I)$  as the initial guess to solve  (\ref{optimalpricing_trip}) using  interior-point algorithm, and obtain the solution $({q}, {\bf {r}}, {\bf {N}}_{A}^I, {\bf {N}}_{H}^I)$  and the corresponding platform profit ${R}$. 
\ENSURE the approximate solution $({q}, {\bf {r}}, {\bf {N}}_{A}^I, {\bf {N}}_{H}^I)$, the corresponding platform profit $R$ and its upper bound $\bar{R}$.
\end{algorithmic}
\end{algorithm}

In the proposed algorithm, we start with an initial guess of the Lagrange multipliers $\mu$, $\kappa$ and solve the globally optimal solution to the sub-problems (\ref{primal_dual_i}), \ref{primal_dual_ii}, and (\ref{primal_dual_q}) in parallel\footnote{Note that both (\ref{primal_dual_i}), (\ref{primal_dual_ii}), and (\ref{primal_dual_q}) are small-scale optimization which can be solved exactly through a grid search.}. The solutions to (\ref{primal_dual_i}), (\ref{primal_dual_ii}), and (\ref{primal_dual_q}) will be used to update $\mu$ based on (\ref{primaldual3}), which will be further utilized to solve (\ref{primal_dual_i}), (\ref{primal_dual_ii}), and (\ref{primal_dual_q}) again. We iterate this process until the convergence criterion is met, and the algorithm is guaranteed to converge because it is a sub-gradient descent algorithm for the dual problem. Depending on the duality gap of (\ref{optimalpricing_trip_relax}), the proposed dual-decomposition algorithm may terminate in two scenarios: (a) there is no duality gap of (\ref{optimalpricing_trip_relax}), and the algorithm converges to a solution that satisfies constraints (\ref{14d}) and (\ref{14f}); (b) there is a duality gap of (\ref{optimalpricing_trip_relax}) and the algorithm terminates to an infeasible point when the algorithm reaches the maximum number of iterations. In either case, the output $\bar{R}$ of our proposed algorithm provides an upper bound on the optimal value of (\ref{optimalpricing_trip}). Part of the solution to (\ref{optimalpricing_trip_relax}), including $\bar{q}, {\bf \bar{r}}$, and ${\bf \bar{N}}_{A}^I$,  will be further used as an initial guess to derive a feasible approximate solution to the original problem (\ref{optimalpricing_trip}) based on a standard interior-point algorithm. At the same time, it offers a theoretical upper bound to quantify the gap between the derived solution (from the interior-point algorithm) to the globally optimal solution to (\ref{optimalpricing_trip}). 

\begin{proposition}
\label{proposition2}
The proposed Algorithm 1 will converge after a finite number of steps. After the algorithm terminates, either (a) $(\bar{q}, {\bf \bar{r}}, {\bf \bar{N}}_{A}^I, {\bf \bar{N}}_{H}^I)$ satisfies constraints (\ref{14d}) and (\ref{14f}), which is a globally optimal solution to (\ref{optimalpricing_trip_relax}), and $\bar{R}$ is an upper bound on the optimal value of (\ref{optimalpricing_trip}); or (b)  $(\bar{q}, {\bf \bar{r}}, {\bf \bar{N}}_{A}^I, {\bf \bar{N}}_{H}^I)$ does not satisfy constraints (\ref{14d}) and (\ref{14f}), but $\bar{R}$ is still an upper bound on the optimal value of (\ref{optimalpricing_trip}).
\end{proposition}
Proposition \ref{proposition2} shows that the proposed algorithm will converge in finite steps and produce an upper bound for the original problem, regardless of whether the relaxed problem is exactly solved or not. The proof of Proposition \ref{proposition2} can be found in \cite[p.385, p.563]{bertsekas1997nonlinear}.  We would like to emphasize that the derived solution $({q}, {\bf {r}}, {\bf {N}}_{A}^I, {\bf {N}}_{H}^I)$ in Algorithm \ref{algorithm1} is a feasible solution to (\ref{optimalpricing_trip}) that satisfies all the constraints (\ref{profit_constraints}). The relaxed problem (\ref{optimalpricing_trip_relax}) is mainly used to provide a good initial guess for solving (\ref{optimalpricing_trip}), and at the same time offer an upper bound on the optimality loss of $({q}, {\bf {r}}, {\bf {N}}_{A}^I, {\bf {N}}_{H}^I)$.

\section{Case Studies}
\label{casestudies}
\begin{figure}[h]
    \centering
    \includegraphics[width=0.65\textwidth]{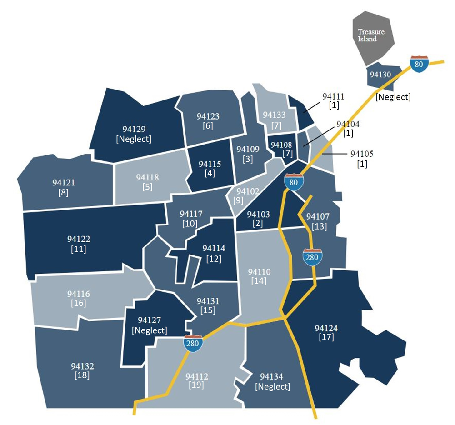}
    \caption{Postal codes and the corresponding zone number in San Francisco. Some zones are neglected in this research due to the negligible trip volume. (Figure courtesy: \url{https://www.usmapguide.com/california/san-francisco-zip-code-map/})}
    \label{zipcode}
\end{figure}

In this section, we conduct simulations based on realistic synthetic data for San Francisco using the aforementioned algorithm. The data set consists of the TNC traffic volume from each origin to each destination in San Francisco at the postal-code level of granularity. It is synthesized based on the total TNC pickup and drop-off data \cite{SFCTA2016} combined with a discrete-choice model calibrated by survey data. A zip-code map for San Francisco is shown in Figure \ref{zipcode}.
Based on the data, we remove zip code zones 94127, 94129, 94130 and 94134  from our analysis as they have negligible trip counts. We further aggregate zip code zones 94111, 94104 and 94105 into a single zone, and aggregate 94133 and 94108 into a single zone, because these zones are close to each other and each of these zones is very small. After this modification, the city of San Francisco is divided into 19 zones. For notation convenience, we number these zones from 1 to 19 and the peer-to-peer matching between the zone number and its zip code is shown in Figure \ref{zipcode} and summarized in \ref{Appendix B}. Note that throughout the case study, we define zone 5, zone 8, zone 11, zone 15, zone 16, zone 17, zone 18, and zone 19 as {\em suburban zones} since they are far from the city center and have fewer passengers. All other zones are  defined  as {\em urban zones}.
{ Furthermore, we assume that the inverse of the average traffic speed\footnote{{ In the interest of computational simplicity, we have elected to employ the inverse of traffic speed (rather than speed itself) in our calculations. This decision stems from the fact that the average trip time $t_{ij}$, as depicted in equation (\ref{traveltime}), is more readily computable when the actual trip time is subject to random fluctuations across diverse trips.}} within the congestion region embodies a linear function of $N_{\mathcal{C}}$, i.e., $\frac{1}{v_c} = \frac{1}{v_c^0}+\rho N_{\mathcal{C}}$, where $v_c^0$ and $\rho$ are considered as model parameters. Under this model, the more vehicles accumulated in the urban areas, the lower the traffic speed, and the longer the travel time. We emphasize that this model is chosen purely for illustrative purposes. The propose modeling framework and solution methodology do not critically depend on the form of the congestion models chosen.}

In summary, the model parameters involved in the numerical study include:
\begin{equation}
    \Theta=\left\{\lambda_{ij}^0,N_0,L,\alpha,\epsilon,\psi,\sigma,\eta,c_{ij}^0,q_0,v_0,\rho,\theta\right\},
\end{equation}
We calibrate these model parameters so that the resultant market outcomes match the real data for San Francisco. Specifically, the potential travel demand is set to satisfy $\lambda_{ij}^0=\tilde{\lambda}_{ij}^H/0.15$, where $\tilde{\lambda}_{ij}^H$ is the observed travel demand in real data corresponding to San Francisco's ride-hailing market with only human drivers. It indicates that about $15\,\%$ of potential passengers choose to take a TNC trip. The traveling time $t_{ij}$ is estimated via Google Maps. The rest of the parameters are set as $N_0 = 10000$, $L=645$, $\alpha=3$, $\epsilon=0.12$, $\sigma=0.17$, $\eta=0.1$, $q_0=\$\,29.34$\,/hour, $\psi = 0.12, v_c^0 = 15$ mph, $\rho=10\%\frac{1}{v_c^0}\frac{1}{1000}$,\footnote{The parameter $\rho$ is chosen so that each 1000 more vehicles will increase the travel time by 10\%.} $v_r = 20$ mph, $\theta = 0.95$.\footnote{It is assumed that the majority of passengers do not exhibit a pronounced inclination towards either AVs or human drivers as providers of ride-hailing services. Only a small percentage of passengers express a preference for being serviced by human drivers.} These values are tuned to render a set of market outcomes that are close to the empirical data of San Francisco (i.e., the ride fare, human driver payment, transport volume, etc.). In particular, after solving (\ref{optimalpricing_trip}) with Algorithm \ref{algorithm1} in the case without AVs (i.e., equilibrium state), we derive that the total passenger arrival rate is about 196 trips/min, the total number of human driver is 3837.54, the driver payment is $\$\,26.55\,$/hour, the average trip fare is about $\$\,44.07$/ride, etc. The passenger demand and driver supply are consistent with results stated in \cite{san2017today}. The driver payment corresponds to the estimation from previous study \cite{parrott2018earnings}. The average trip fare matches the report from Lyft in 2021 \cite{Lyft}. Based on these parameters, we execute the proposed Algorithm \ref{algorithm1} in $Matlab^@$ 2022a on a Dell desktop with 8-core i7-9700 CPU. Each execution of Algorithm  \ref{algorithm1} takes about 172.45 seconds. The following subsections present three rounds of simulations that validate the proposed  framework: one round without any regulation (Section \ref{subsection5_1}), the second round under minimum wage regulations (Section \ref{subsection5_2}), and the third round under pickup restrictions on AVs (Section \ref{sectin 5.3}).

\subsection{Evaluation of Transport Equity}
In the case study, we will examine the impacts of AVs on transport equity. To acquaint readers with the concept utilized in the simulation study, we will provide a rigorous definition of transport equity in this subsection. In particular, equity concerns the distribution of impacts among populations \cite{litman2002evaluating}. It can be examined from {\em horizontal} and {\em vertical} perspectives. Horizontal equity requires that individuals with comparable abilities or needs be treated equally in the distribution of benefits and costs, e.g., people in different regions of the city should have access to the same quality of mobility services (this is also referred to as spatial equity). On the other hand, vertical equity requires that the allocation of benefits and costs favors disadvantaged groups, e.g., people that are disadvantaged in social-economic status (e.g., elderly people, the disabled, etc) should be favored by transport policies. This is also referred to as social equity.    To quantify transport equity, the first step involves segmenting the population into subgroups based on geographic location (for horizontal equity)  or socio-economic status (for vertical equity). After the segmentation, the costs and benefits of each group will be identified \cite{guo2020systematic},  which are then used to quantify transport equity across distinct populations using inequity indicators (e.g., Gini index, Theil index, Atkinson index, etc).

We use Theil coefficient \cite{theil1967economics} to quantify the transport equity in the multimodal transport system. To this end, we first define $u_{ij}^{AH}=-c_{ij}^{AH}$ and  $u_{ij}^{H}=-c_{ij}^{H}$ as the average utility of passengers from origin zone $i$ to destination zone $j$ for Class 1 and Class 2 passengers, respectively. Denote $u_{ij}^{0}=-c_{ij}^{0}$ as the average utility of the outside options.  Under these notations, we further define accessibility measure $A_{ij}^{AH}$ and $A_{ij}^{H}$ as the expected maximum utility for Class 1 and Class 2 passengers from origin $i$ to destination $j$ in the multimodal transportation system, respectively. Based on the logit model (\ref{2}), the maximum expected utility $A_{ij}^{AH}$ and $A_{ij}^{H}$ can be written as the logsum formula \cite{ho2006combined,van2021evaluating}:
\begin{equation}
    \begin{cases}
         A_{ij}^{AH} = \frac{1}{\epsilon} \log \left( \exp(\epsilon u_{ij}^{AH})+\exp(\epsilon u_{ij}^0) \right), \\
         A_{ij}^{H} = \frac{1}{\psi} \log \left( \exp(\psi u_{ij}^{H})+\exp(\psi u_{ij}^0) \right).
    \end{cases}
\label{logsum_accessibility}
\end{equation}
Given the accessibility measures, we further denote $\lambda_{i}^{AH}$ and $\lambda_{i}^{H}$ as the arrival rates of potential demand of Class 1 and Class 2 passengers from zone $i$, respectively; let  $\lambda^{AH}$ and $\lambda^{H}$ be arrival rate of potential demand of passengers of Class 1 and Class 2, respectively, regardless of where they start their trip; let $\overline{\lambda}$ denote the total passenger demand. They satisfy:
\begin{equation}
\begin{cases}
    \lambda_{i}^{AH} = \sum_{j=1}^M \theta \lambda_{ij}^0; \\
    \lambda_{i}^{H} = \sum_{j=1}^M (1-\theta) \lambda_{ij}^0;\\
    \lambda^{AH} = \sum_{i=1}^M \sum_{j=1}^M \theta \lambda_{ij}^0; \\
    \lambda^{H} = \sum_{i=1}^M \sum_{j=1}^M (1-\theta) \lambda_{ij}^0;  \\
    \overline{\lambda} = \sum_{i=1}^M \sum_{j=1}^M  \lambda_{ij}^0.
\end{cases}\end{equation}Correspondingly, we denote $A_{i}^{AH}$ and $A_{i}^{H}$ as the average accessibility of passengers originating from zone $i$ in Class 1 and Class 2, respectively. They can be calculated as:
\begin{equation}
\begin{cases}
    A_{i}^{AH} = \frac{\sum_{j=1}^M \left(   \theta \lambda_{ij}^0 A_{ij}^{AH} \right)}{\lambda_{i}^{AH}}; \\
    A_{i}^{H} = \frac{\sum_{j=1}^M   \left( (1-\theta)\lambda_{ij}^0 A_{ij}^H\right)}{\lambda_{i}^H}. 
\end{cases} \end{equation}Let $A^{AH}$ and $A^H$ denote the average accessibility of passengers in Class 1 and Class 2, respectively, regardless of where they start their trips. They can be calculated as:
\begin{equation}
\begin{cases}
    A^{AH} = \frac{\sum_{i=1}^M\sum_{j=1}^M \left(   \theta \lambda_{ij}^0 A_{ij}^{AH} \right)}{\lambda^{AH}};\\
    A^H = \frac{\sum_{i=1}^M \sum_{j=1}^M   \left( (1-\theta)\lambda_{ij}^0 A_{ij}^H\right)}{\lambda^{H}}. 
\end{cases}
\end{equation}Furthermore, let $\bar{A}$ denote the average accessibility of all passengers. It can be calculated as:
\begin{equation}
\overline{A} = \frac{\sum_{i=1}^M\sum_{j=1}^M \left(   \theta \lambda_{ij}^0 A_{ij}^{AH} \right)+\sum_{i=1}^M \sum_{j=1}^M   \left( (1-\theta)\lambda_{ij}^0 A_{ij}^H\right) }{\bar{\lambda}},\end{equation}Finally, we define the Theil (T) coefficient of the accessibility distribution in the multi-modal transportation system as:
\begin{align} \label{Theil_coefficient}
    T = &\underbrace{\sum_{i=1}^M \left(\frac{\lambda_{i}^{AH}}{\overline{\lambda}} \frac{A_{i}^{AH}}{A^{AH}}\right) \ln \left(\frac{A_{i}^{AH}}{A^{AH}}\right)   +  \sum_{i=1}^M \left(\frac{\lambda_{i}^{H}}{\overline{\lambda}} \frac{A_{i}^{H}}{A^{H}}\right) \ln \left(\frac{A_{i}^{H}}{A^{H}}\right)    }_\text{Spatial equity}  \\ \notag
    &+  \underbrace{ \left(\frac{\lambda^{AH}}{\bar{\lambda}} \frac{A^{AH}}{\overline{A}}\right) \ln \left(\frac{A^{AH}}{\bar{A}}\right)   
     +  \left(\frac{\lambda^{H}}{\bar{\lambda}} \frac{A^{H}}{\overline{A}}\right) \ln \left(\frac{A^{H}}{\bar{A}}\right)}_\text{Social equity}.
\end{align}

As shown in (\ref{Theil_coefficient}), the Theil coefficient can be decomposed into two components. The first two terms are also referred to as the ``WITHIN'' component, which calculates the inequality in the distribution of passenger accessibility across distinct geographic zones {\em within} each group. Consequently, it characterizes the spatial inequity due to the differentiated mobility service across distinct zones. The last two terms are also referred to as the  ``BETWEEN'' component evaluates the difference in the distribution of average passenger accessibility {\em between} different classes, which captures the social inequity arising from the differences of passengers in their status (e.g., regular, elderly, disabled, etc). Note that a larger Theil coefficient indicates a more inequitable distribution of benefits across different zones or socioeconomic groups. It can be easily verified that if all zones have the same level of accessibility, then the ``WITHIN'' component of (\ref{Theil_coefficient}) is zero. Similarly, if all classes of passengers have the same level of accessibility, then the ``BETWEEN'' component of (\ref{Theil_coefficient})  is zero.

\subsection {TNC Market in the Absence of Regulations}
\label{subsection5_1}
\begin{figure}[h!]
    \centering
    \begin{minipage}[h]{0.49\textwidth}
    \centering
    \includegraphics[width=1\textwidth]{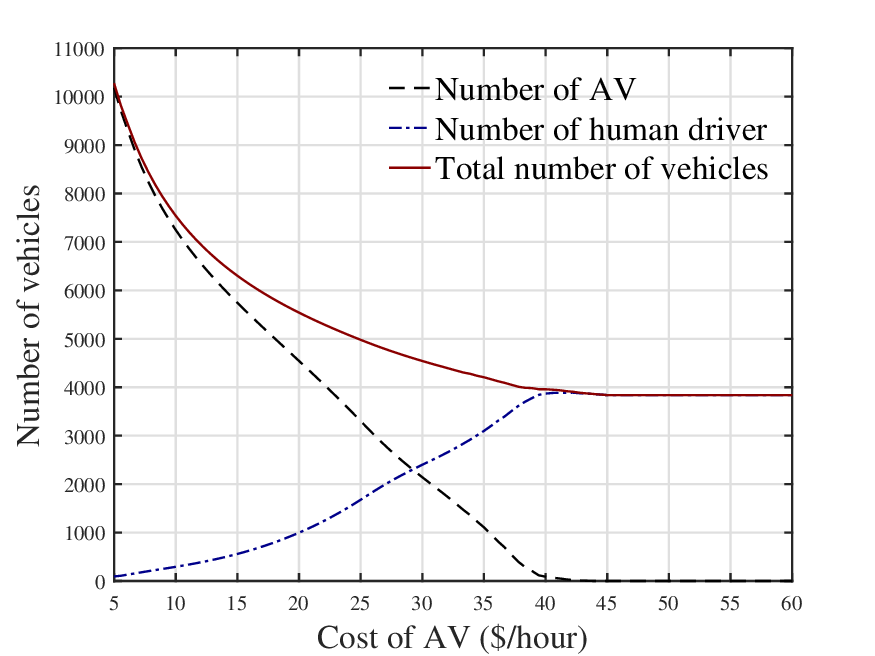}
    \caption{Number of TNC vehicles under different values of AV cost}
    \label{totalnovehicle}   
    \end{minipage}
    \begin{minipage}[h]{0.49\textwidth}
    \centering
    \includegraphics[width=1\textwidth]{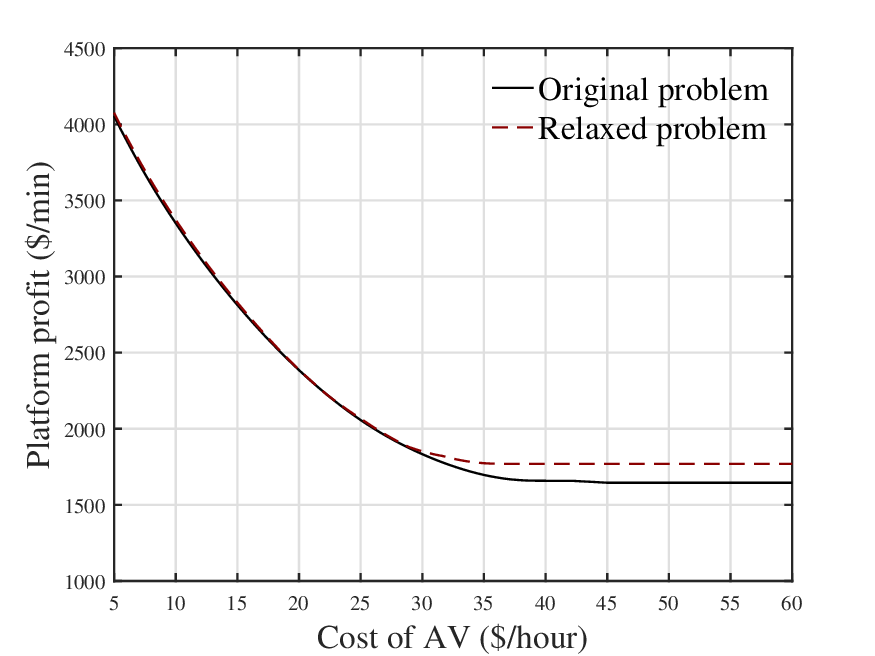}
    \caption{Comparison of TNC profit for the original and relaxed problem under varying AV cost}
    \label{Comparison_relax_r1}
    \end{minipage}
\end{figure}

\begin{figure}[!ht]
    \centering
    \vspace{0cm}
    \begin{minipage}[h]{0.49\textwidth} 
    \centering
    \includegraphics[width=1\textwidth]{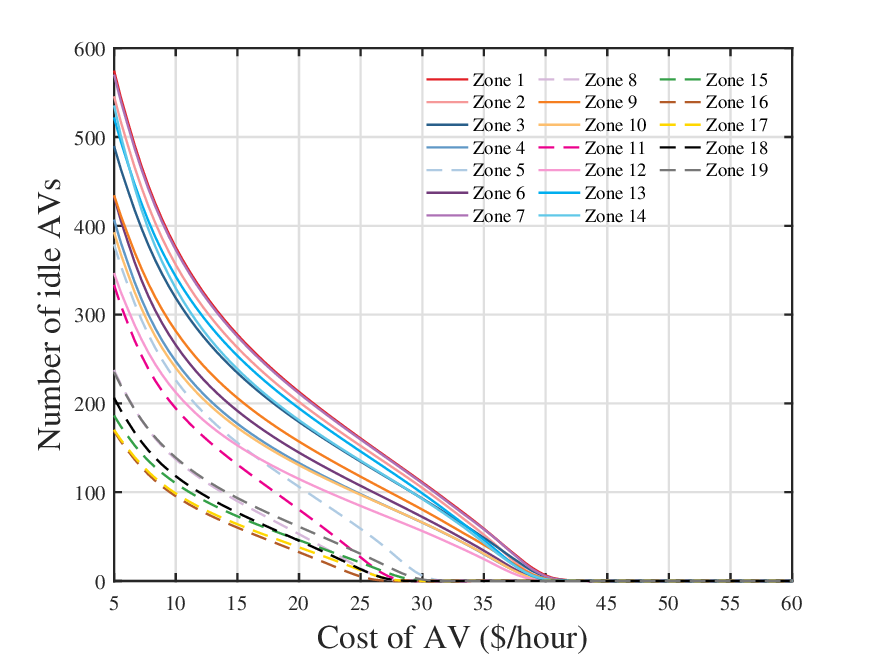}
    \caption{Number of idle AVs in each zone under varying cost of AV}
    \label{idleAV_figure2}
    \end{minipage}
    \begin{minipage}[h]{0.49\textwidth} 
    \includegraphics[width=1\textwidth]{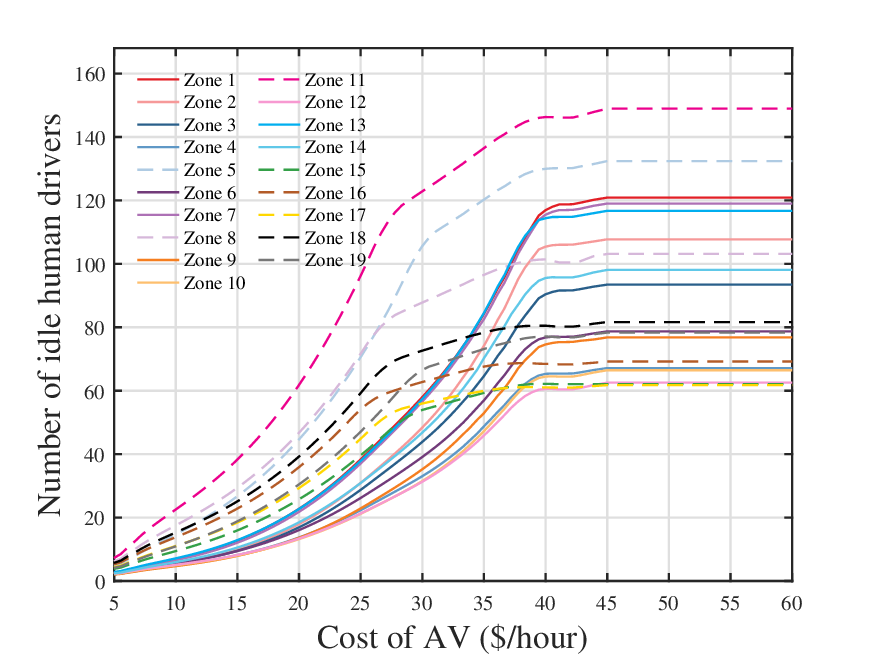}
    \caption{Number of idle human drivers in each zone under varying cost of AV}
    \label{idleHDV}
    \end{minipage}
    \begin{minipage}[h]{0.49\textwidth} 
    \includegraphics[height=6.4cm,width=8.7cm]{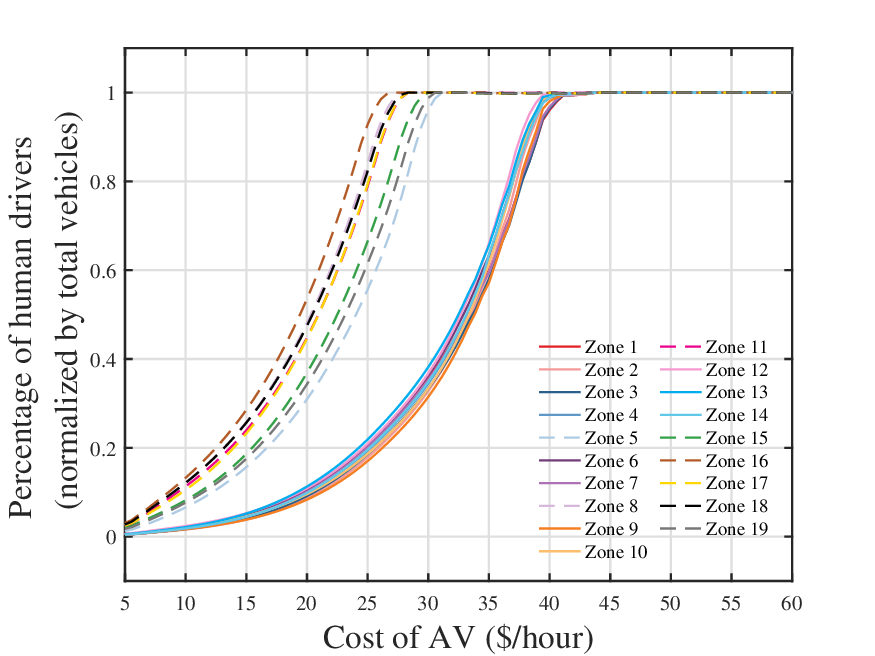}
    \caption{Percentage of human drivers in each zone under varying cost of AV (normalized by total vehicles)}
    \label{percentage of HD}
    \end{minipage}
    \begin{minipage}[h]{0.49\textwidth} 
    \includegraphics[height=6.4cm,width=8.7cm]{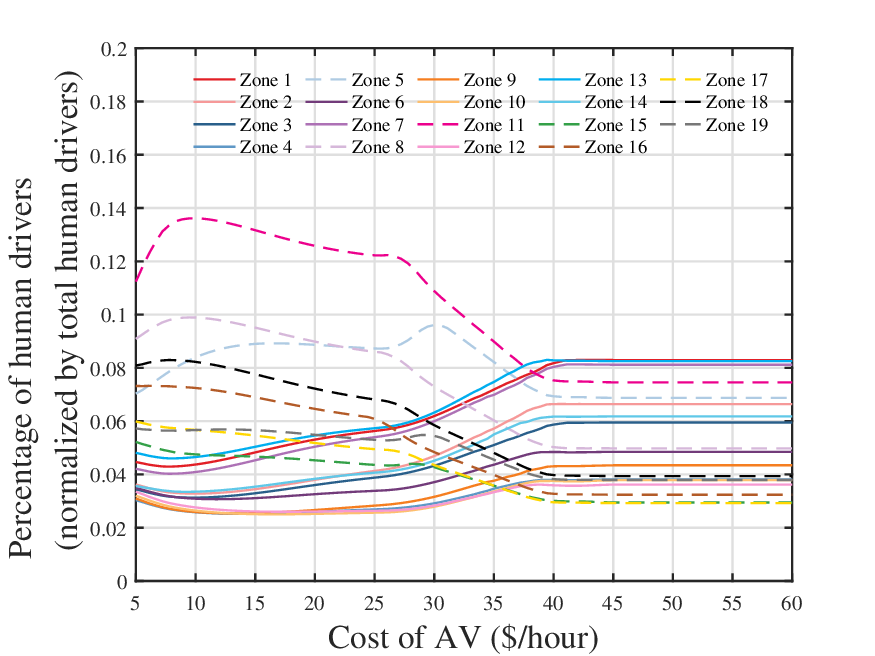}
    \caption{Percentage of human drivers in each zone under varying cost of AV (normalized by total human drivers)}
    \label{percentage of HD2 r1}
    \end{minipage}
    \begin{minipage}[h]{0.49\textwidth} 
    \includegraphics[height=6.4cm,width=8.7cm]{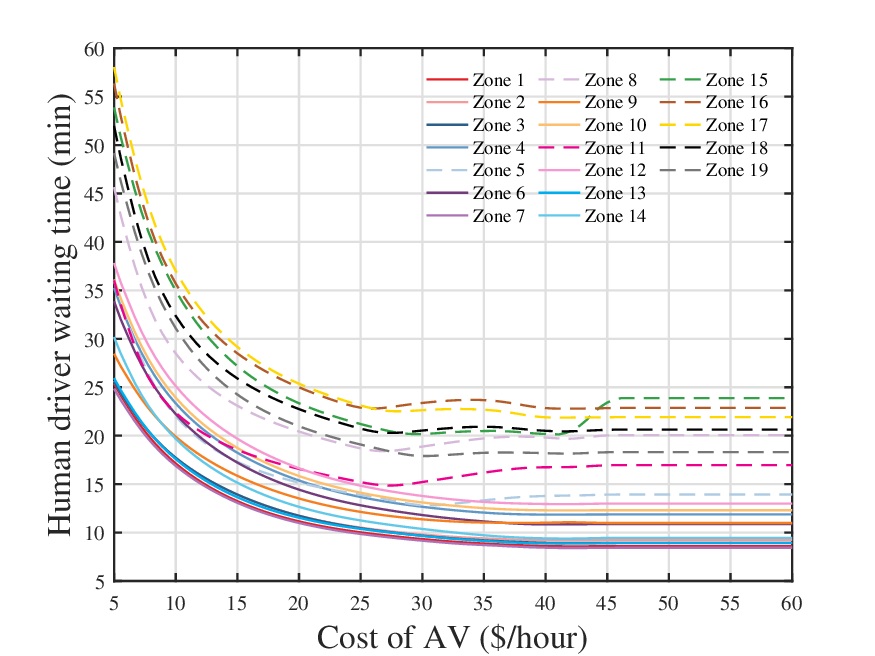}
    \caption{Human driver waiting time in each zone under varying cost of AV}
    \label{driverwait}
    \end{minipage}
    \begin{minipage}[h]{0.49\textwidth} 
    \includegraphics[height=6.4cm,width=8.7cm]{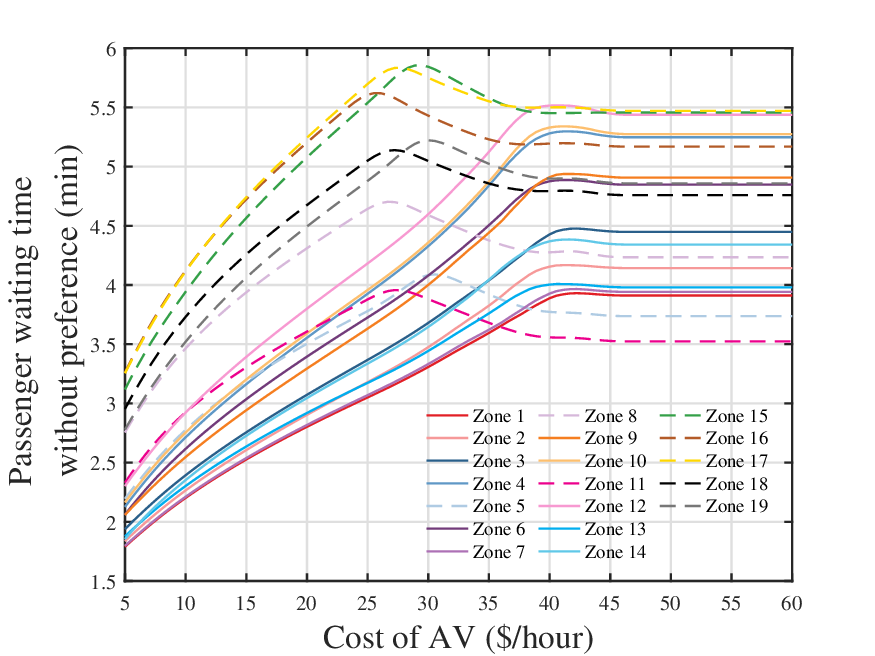}
    \caption{Passenger waiting time (without preference) in each zone under varying cost of AV}
    \label{passengerwait}
    \end{minipage}
    \end{figure}

\begin{figure}
    \begin{minipage}[h]{0.49\textwidth} 
    \includegraphics[height=6.4cm,width=8.7cm]{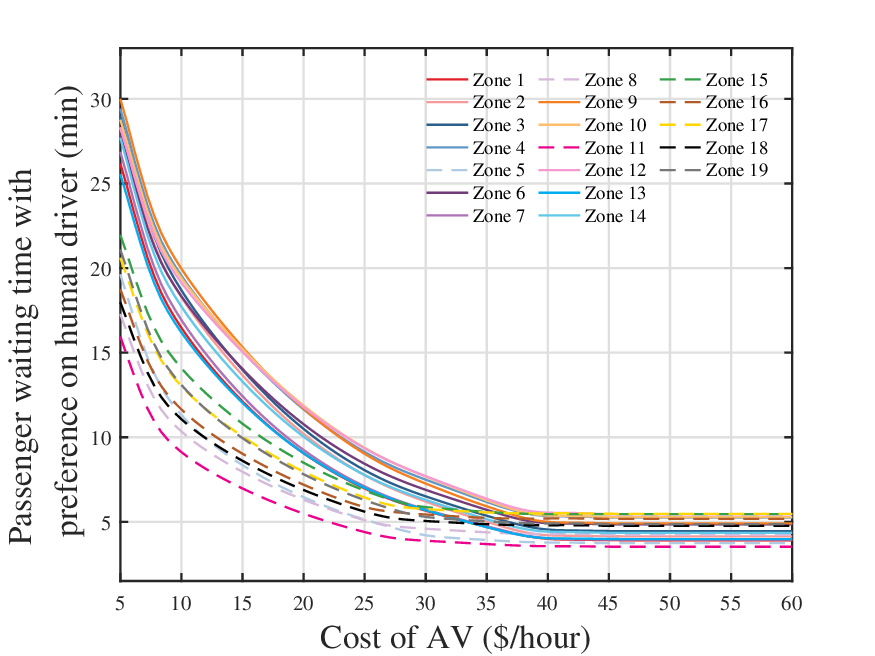}
    \caption{Passenger waiting time (with preference on human driver) in each zone under varying cost of AV}
    \label{passengerwait1}
    \end{minipage}
    \begin{minipage}[!ht]{0.49\textwidth} 
    \includegraphics[height=6.4cm,width=8.7cm]{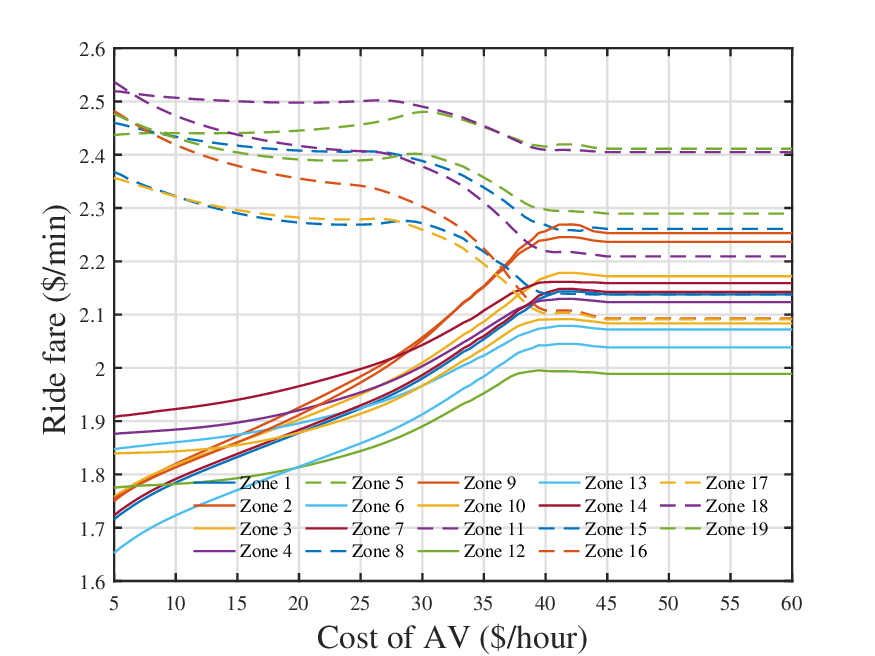}
    \caption{Ride fares per-minute in each zone under varying cost of AV}
    \label{ridefare}
    \end{minipage}
\end{figure}

In this subsection, we solve the platform's profit maximization problem (\ref{optimalpricing_trip}) in the absence of regulations.  We fix all other model parameters and vary the cost of AV, i.e., $C$, to investigate how the AV deployment cost will affect the market outcomes in the mixed environment. As the outside option of human drivers corresponds to a payment $q_0=\$\,29.34$/hour, we perturb $C$ between $\$\,5/$hour and $\$\,60/$hour to cover a wide range of AV costs. The proposed Algorithm \ref{algorithm1} was executed and the computed solution $(\bar{q}, {\bf \bar{r}}, {\bf \bar{N}}_{A}^I, {\bf \bar{N}}_{H}^I)$ satisfies constraints (\ref{14d}) and (\ref{14e}). Using Proposition \ref{proposition2}, we find that the obtained solution $(\bar{q}, {\bf \bar{r}}, {\bf \bar{N}}_{A}^I, {\bf \bar{N}}_{H}^I)$ is the global optimal solution to problem (\ref{optimalpricing_trip_relax}), which provides a performance guarantee for the approximate solutions $({q}, {\bf {r}}, {\bf {N}}_{A}^I, {\bf {N}}_{H}^I)$. Figure \ref{totalnovehicle}-\ref{ridefare} present the market outcomes (including the vehicle supply, driver wage, ride fare, platform profit, etc.) as a function of $C$.

Figure \ref{totalnovehicle}-\ref{Comparison_relax_r1} show the aggregate market outcomes under different values of AV cost. Based on Figure~\ref{totalnovehicle},  AVs and human drivers will co-exist in the TNC fleet when $\$\,5\text{/hour} < C < \$\,42.5\text{/hour}$. Based on Figure~\ref{Comparison_relax_r1}, we can evaluate the gap between the derived solution and the unknown globally optimal solution. It is clear that the platform profit derived by the proposed algorithm is very close to that of the globally optimal solution to the non-concave program (\ref{optimalpricing_trip}): when the cost of AV is large, the gap between these platform profits at two solutions is around $7.01\,\%$; when the cost of AV is small, the gap between these two cases is almost negligible. Note that this is because when the cost is small, there is a large number of AVs in the mixed fleet. In this case, the platform can satisfy the rebalancing constraints for human drivers (\ref{10f}) and AVs (\ref{10g}) at a negligible cost by jointly determining the spatial distribution of idle AVs and the repositioning flow of AVs. For instance, for a given human-driver repositioning flow $\Tilde{f}_{ij,H}$, the platform can first adjust $N_{i,A}^I$ for each zone so that (\ref{10f}) is satisfied, and then adjust $f_{ij,A}$ for each pair of $i$ and $j$ so that (\ref{10g}) is satisfied. As the number of controllable variables (i.e., $N_{i,A}^I$ and $f_{ij,A}$) is large, these adjustments can be performed with little impact on the platform profit.  However, when the cost of AVs is large, there is a small number of AVs (or no AVs) in the fleet, and it is difficult to satisfy (\ref{10f}) and (\ref{10g}) by adjusting the value of $N_{i,A}^I$ and $f_{ij,A}$. This leads to a larger gap between the optimal value of (\ref{optimalpricing_trip}) and that of (\ref{optimalpricing_trip_relax}).

\begin{figure}[t!]
    \centering
    \begin{minipage}[h]{0.49\textwidth}
    \centering
    \includegraphics[width=1\textwidth]{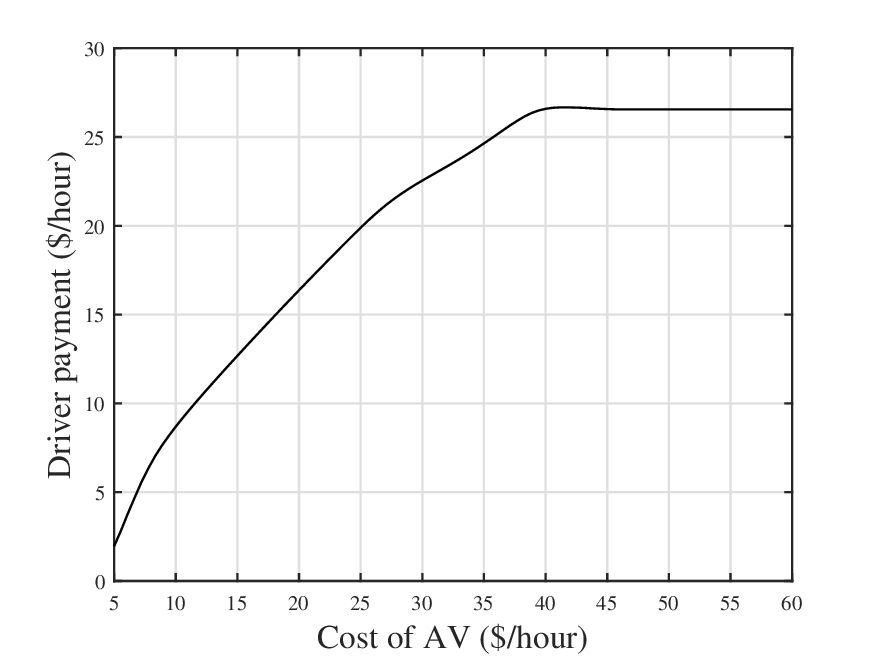}
    \caption{Hourly wage of human drivers under different values of AV cost}
    \label{driverpaymentr1}   
    \end{minipage}
    \begin{minipage}[h]{0.49\textwidth}
    \centering
    \includegraphics[width=1\textwidth]{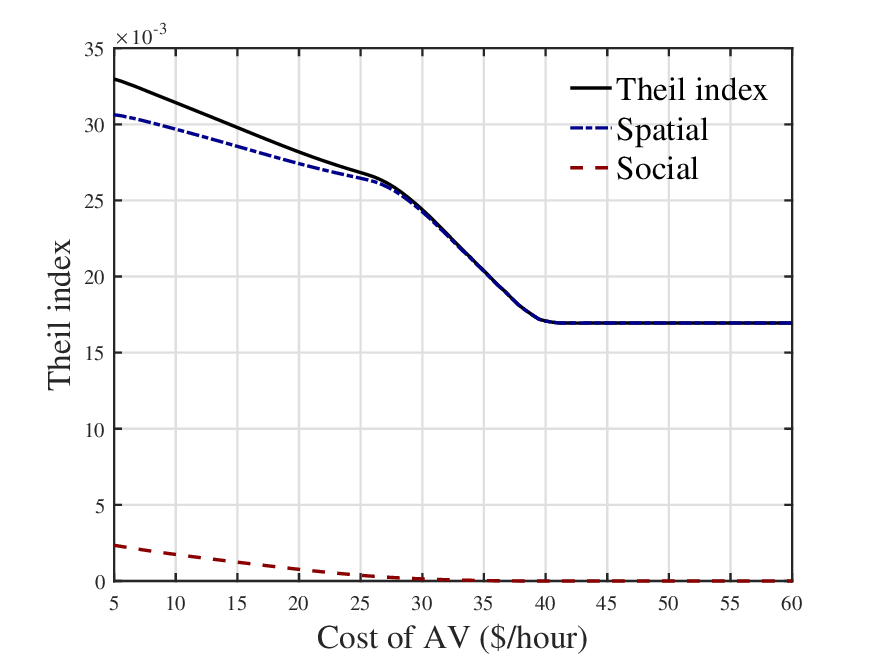}
    \caption{Theil coefficient, Spatial equity, and Social equity under varying cost of AV}
    \label{covr1}
    \end{minipage}
\end{figure}

Figure \ref{idleAV_figure2}-\ref{ridefare} show the market outcomes in different zones under different values of $C$. As mentioned before, we divide the city of  San Francisco into 19 zones and categorize each zone as urban zone or suburban zone based on its proximity to the center of the city\footnote{Suburban zones include zone 5, zone 8, zone 11, zone 15, zone 16, zone 17, zone 18, and zone 19.}. The results of suburban zones are presented with dashed lines, while the results of the urban zones are presented with solid lines. Interestingly, we observe the following results in the simulation studies:
\begin{itemize}
    \item As the cost of AV reduces, the platform activates AVs in different zones in different orders (Figure \ref{idleAV_figure2}). When $C\geq \$\,42.5$/hour, the platform does not deploy any AVs in the TNC fleet. As $C$ reduces (horizontal axis of Figure \ref{idleAV_figure2} changes from right to left), the platform first deploys AVs in zone 1 at around $C=\$\,42.5$/hour (red solid curve in Figure \ref{idleAV_figure2}), last deploys AVs in zone 16 at around $C=\$\,25.5$/hour (yellow dashed curve in Figure \ref{idleAV_figure2}), and distributes AVs in other zones sequentially in between, i.e., $\$\,5\text{/hour }\leq C\leq \$\,42.5\text{/hour}$. By comparing the dashed lines with solid lines in Figure \ref{idleAV_figure2}, it is clear that the platform has a tendency to prioritize AV deployment in urban areas in stead of suburban areas. This is because urban areas have higher passenger demand, and deploying AVs in these areas will bring a higher profit to the platform.
    
  
    \item As AVs flourish in urban areas, they will replace the majority of human drivers and force them to relocate to other zones. Therefore, the spatial distribution of AVs and human drivers are highly asymmetric at the profit-maximizing equilibrium. On the one hand, as shown in Figure \ref{idleAV_figure2}, AVs will be highly concentrated in the urban areas where the passenger demand is higher. On the other hand, as shown in Figures \ref{idleHDV}, \ref{percentage of HD} and \ref{percentage of HD2 r1}, it is clear that in these figures, all dashed lines are above the solid lines when AV cost is small, thus we can conclude the AVs will push human drivers to relocate to suburban areas where the earning opportunity is much lower when the cost of AV is small.  Figure \ref{driverwait} presents the driver waiting time in each zone. This reflects that human drivers are disproportionately concentrated in suburban areas and they are forced to take a higher waiting time in these areas, which is unfair.
    \item The introduction of AVs will exacerbate the spatial inequity of service qualities among passengers. Figure \ref{passengerwait} shows the average passenger waiting time in each zone. For the majority of human passengers who do not have a preference between human drivers and AVs, as the cost of AVs decreases, passengers benefit from shorter waiting times. However, it should be noted that the waiting time in remote areas is generally higher compared to that in urban areas, especially after AVs are introduced: when the cost parameter $C$ is smaller than \$30 per hour, the majority of the dashed lines in the graph are  positioned above the solid lines, indicating a disparity in waiting times between urban and suburban areas following the introduction of AVs. Conversely, for the small fraction of human passengers who have a preference for human drivers, their waiting time is significantly increased as AV adoption increases (see Figure \ref{passengerwait1}). Figure \ref{ridefare} shows the average passenger ride fare in each zone of the city. When $C$ is small, we can again observe that all dashed lines are above the solid lines, and the ride fare may even increase as $C$ reduces. This verifies a significant gap of ride fares between urban and suburban areas when the cost of AV is smaller than 35\$ per hour. Overall, Figure \ref{passengerwait}-\ref{ridefare} reveal that the deployment of AVs may intensify the inequity gap of travel cost between distinct areas and between distinct classes of passengers.
\end{itemize}

Based on the above discussion, it is clear that deployment of AVs  will negatively affect human drivers (see Figure \ref{driverpaymentr1}) and exacerbate the existing equity gaps among TNC passengers. To make this argument more precise, in this numerical study, we calculate the Theil index under distinct value of AV cost, and show how the distribution of benefits among passengers will changes as AVs adoption increases. The results are shown in Figure \ref{covr1}, where both the ``WITHIN'' component, ``BETWEEN'' component, and the overall index are shown. It is clear that as the cost of AV reduces from $\$\,60$/hour to $\$\,5$/hour, the Theil index increases from 0.0169 to 0.033 , with the ``WITHIN'' component increasing from 0.0169 to 0.0306, and the ``BETWEEN'' component increasing from 0 to 0.0024 This clearly indicates that increasing the scale of AVs in the mixed fleet will lead to exacerbated transport inequity, both from the spatial perspective (``WITHIN'' component) and the social perspective (``BETWEEN'' component).

The above observations motivate us to consider regulatory policies that protect human drivers and improve transport equity. To this end, we propose two relevant policies: (a) a minimum driver wage on human drivers, and (b) a restrictive pickup policy that prohibits AVs to pick up passengers from urban areas\footnote {Note that they are permitted to drop off passenger in urban areas.}, which is consistent with New York City's policy that prevents green taxis from picking up passengers in South Manhattan. The rest of this section evaluates the impacts of these two regulatory policies.

\subsection{Impacts of Minimum Wage for Human Drivers}
\label{subsection5_2}

To protect human drivers from the negative impacts of AVs, we consider a wage floor on human driver earnings. This can be captured as an additional constraint, i.e., $q\geq q_{min}$, where $q_{min}$ is the minimum wage required by the government.  Under the wage floor, the profit-maximizing problem of the platform can be written as:
\begin{equation}
\hspace{-5cm}
\label{optimalpricing_trip_regulated}
    \mathop {\max }\limits_{ q, {\bf r}, {\bf N}_{A}^I, {\bf f}_A }\sum_{i=1}^{M}\sum_{j=1}^{M}{r_i\left(\lambda_{ij}^{AH}+\lambda_{ij}^H\right)t_{ij}}-C\cdot N_{A}-N_0F_d\left(q\right)q.
\end{equation}
\begin{subnumcases}
   {\label{profit_constraints_regulated}}
   \lambda_{ij}^{AH} =\theta\lambda_{ij}^0  \frac{e^{-\epsilon c_{ij}^{AH}}}{e^{-\epsilon c_{ij}^{AH}}+e^{-\epsilon c_{ij}^0}},
    \label{10a_reg}\\
    \lambda_{ij}^{H} = (1-\theta)\lambda_{ij}^0  \frac{e^{-\psi c_{ij}^{H}}}{e^{-\psi c_{ij}^{H}}+e^{-\psi c_{ij}^0}},
    \label{10aa_reg}\\
     w_i^{p,AH}(N_i^I) \leq w_i^{p,max}, \quad w_i^{p,A}(N_i^{I,A}) \leq w_i^{p,max},
    \label{10b_reg}\\
      N_{i,A}^I = w_i^{d,A}\cdot \displaystyle \sum_{j=1}^{M}\left(\dfrac{N_{i,A}^I}{N_{i,A}^I+N_{i,H}^I} \lambda_{ij}^{AH} \right), 
    \label{10d_Reg}\\
     N_{i,H}^I = w_i^{d,H}\cdot \displaystyle\sum_{j=1}^{M}\left( \dfrac{N_{i,H}^I}{N_{i,A}^I+N_{i,H}^I} \lambda_{ij}^{AH} + \lambda_{ij}^{H}\right), \label{10dd_Reg}\\
     \sum_{j=1}^{M}\left({\dfrac{N_{j,A}^I}{N_{j,A}^I+N_{j,H}^I}}\lambda_{ij}^{AH}+\tilde{f}_{ji,A}\right)=\sum_{j=1}^{M}\left(\dfrac{N_{i,A}^I}{N_{i,A}^I+N_{i,H}^I}\lambda_{ij}+\tilde{f}_{ij,A}\right),
    \label{10f_Reg}\\
    \sum_{j=1}^{M}\left(\dfrac{N_{j,H}^I}{N_{j,A}^I+N_{j,H}^I}\lambda_{ji}^{AH}+\lambda_{ji}^{H}+\Tilde{f}_{ji,H}\right)=\sum_{j=1}^{M}\left(\dfrac{N_{i,H}^I}{N_{i,A}^I+N_{i,H}^I}\lambda_{ij}^{AH}+\lambda_{ij}+\Tilde{f}_{ij,H}\right),
    \label{10g_Reg}\\
    {N_0}{F_d}(q)\geq \displaystyle\sum_{i=1}^{M}{\sum_{j=1}^M \dfrac{N_{i,H}^I}{N_{i,A}^I+N_{i,H}^I}\lambda_{ij}^{AH}t_{ij}}+\displaystyle\sum_{i=1}^M\sum_{j=1}^M \lambda_{ij}^{H}t_{ij} \\\notag
      \quad \quad  \quad \quad \quad  +\displaystyle\sum_{i=1}^{M}\sum_{j=1}^{M} {\dfrac{N_{i,H}^I}{N_{i,A}^I+N_{i,H}^I}\lambda_{ij}^{AH}w_i^{p,AH}}+\sum_{i=1}^{M}\sum_{j=1}^{M}\lambda_{ij}^{H}w_i^{p,H}+\sum_{i=1}^{M}N_{i,H}^I,
    \label{10h_reg}\\
    q\geq q_{min}, \label{10i_reg}\\
    N_{A}=\displaystyle\sum_{i=1}^M\sum_{j=1}^M \dfrac{N_{i,A}^I}{N_{i,A}^I+N_{i,H}^I}\lambda_{ij}^{AH}t_{ij} +\displaystyle\sum_{i=1}^{M}\sum_{j=1}^{M} \dfrac{N_{i,A}^I}{N_{i,A}^I+N_{i,H}^I}\lambda_{ij}^{AH}w_i^{p,AH}
    +\sum_{i=1}^{M}N_{i,A}^I,
    \label{10j_Reg}\\
    N_\mathcal{C} = \displaystyle\sum_{i=1}^{M}{\sum_{j=1}^M\lambda_{ij}}d_{ij}^\mathcal{C}\cdot \frac{1}{v_c}+\sum_{i\in \mathcal{C}}\sum_{j=1}^{M}{\lambda_{ij}w_i^p}+\sum_{i\in \mathcal{C}}{\left(N_{i,A}^I+N_{i,H}^I\right).}
    \label{10k_Reg}
\end{subnumcases}

Compared to (\ref{optimalpricing_trip}), the major difference is that we impose the minimum wage constraint (\ref{10i_reg}) and relax the equality constraint (\ref{10h}) to an inequality constraint (\ref{10h_reg}), which allows the platform to choose a subset of drivers who are willing to enter the industry under the minimum wage regulation\footnote{The left-hand side of (\ref{10h_reg}) is the number of drivers who are willing to work for the platform, and the right-hand side is the number of drivers that are actually hired by the platform. Under a minimum wage, the platform may find it too expensive to hire all the willing drivers. In this case, it may only hire a subset of drivers, leading to the inequality in (\ref{10h_reg}).}. 

\begin{figure}[h!]
    \begin{minipage}[t]{0.48\textwidth}
    \centering
    \includegraphics[width=1\textwidth]{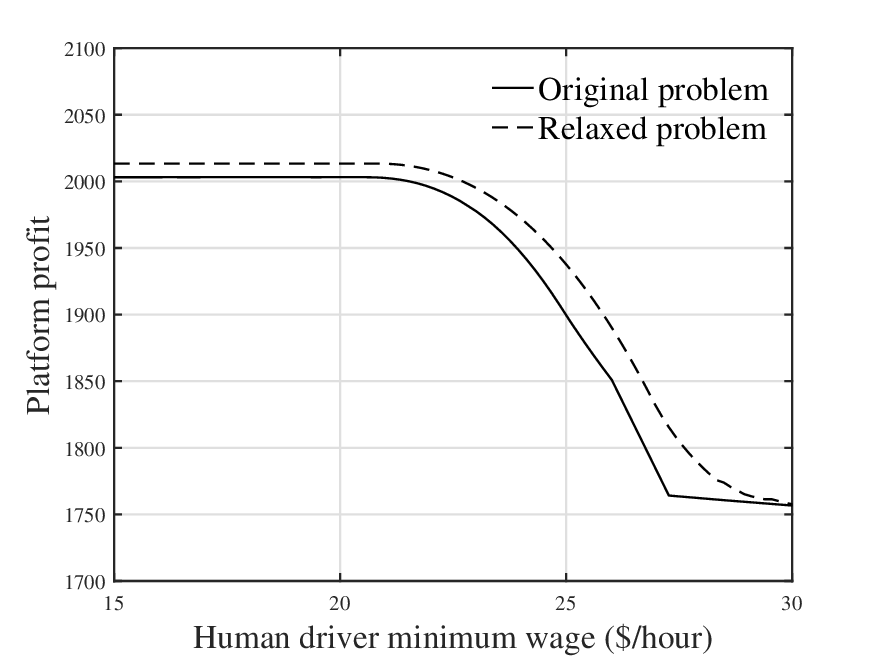}
    \caption{Comparison of TNC profit for original and relaxed problems under varying $q_{min}$}
    \label{Comparison_relax_r2}
    \end{minipage}
    \begin{minipage}[t]{0.48\textwidth} 
    \centering
    \includegraphics[width=1\textwidth]{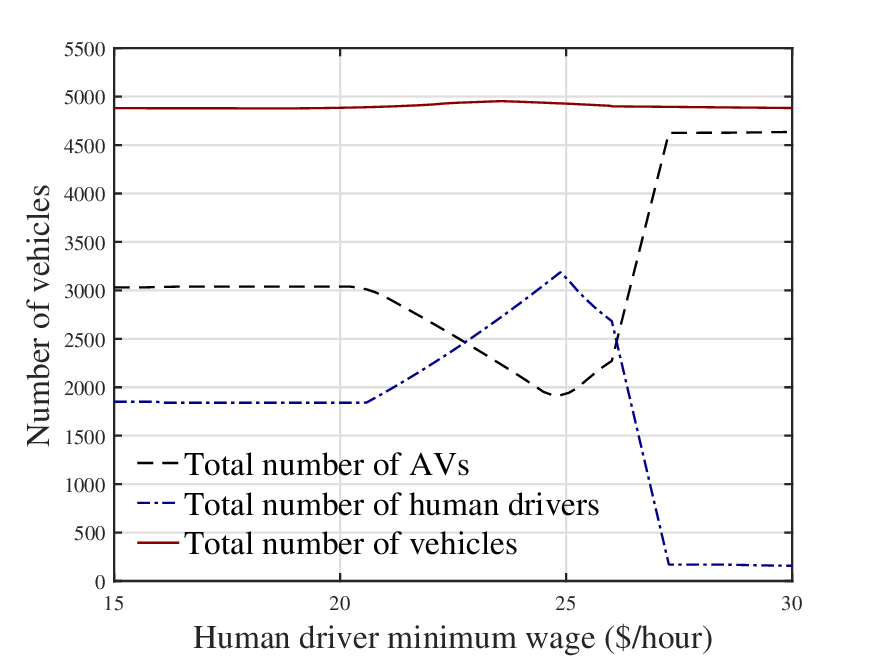}
    \caption{The number of vehicles under varying minimum driver payment $\it{q_{min}}$}
    \label{totalnovehicle_r2}
    \end{minipage}
    \begin{minipage}[!htbp]{0.49\textwidth}
    \includegraphics[height=6.4cm,width=8.7cm]{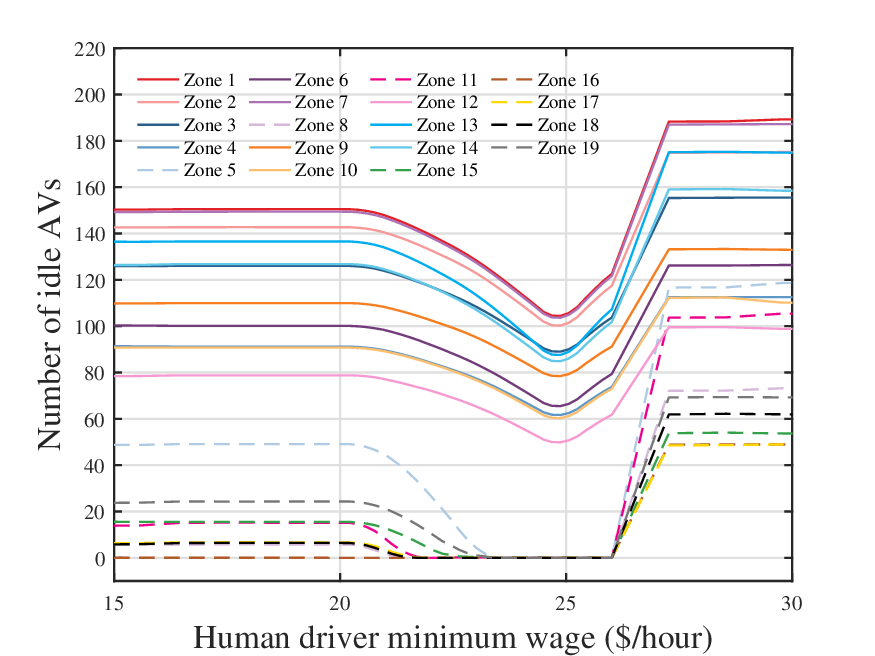}
    \caption{Number of idle AV in each zone under varying minimum driver payment $\it{q_{min}}$}
    \label{idleAV_r2}
    \end{minipage}
    \begin{minipage}[!htbp]{0.49\textwidth}
    \includegraphics[height=6.4cm,width=8.7cm]{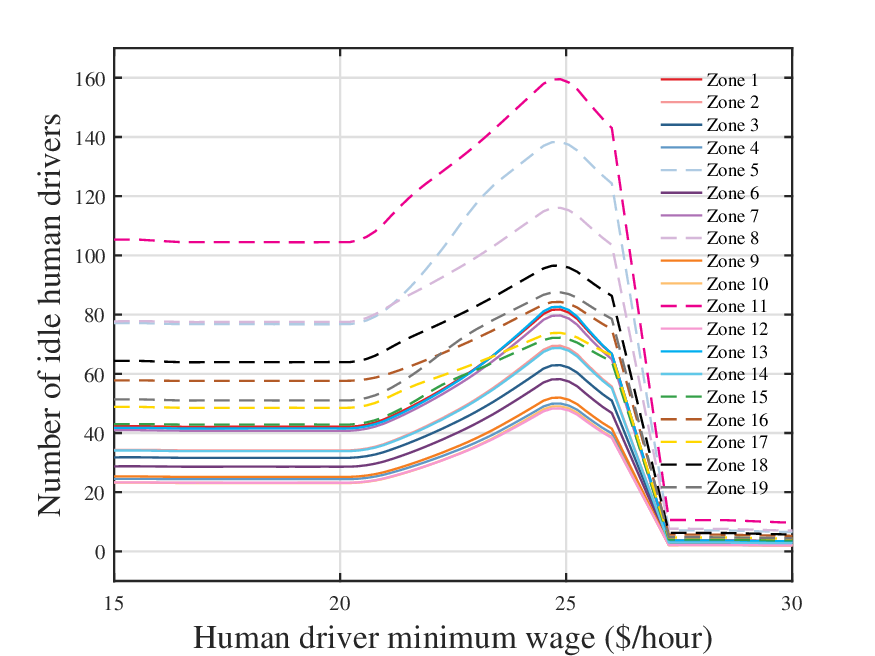}
    \caption{Number of idle human driver in each zone under varying minimum driver payment $\it{q_{min}}$}
    \label{idleHDV_r2}
    \end{minipage}
    \begin{minipage}[!htbp]{0.49\textwidth}
    \includegraphics[height=6.4cm,width=8.7cm]{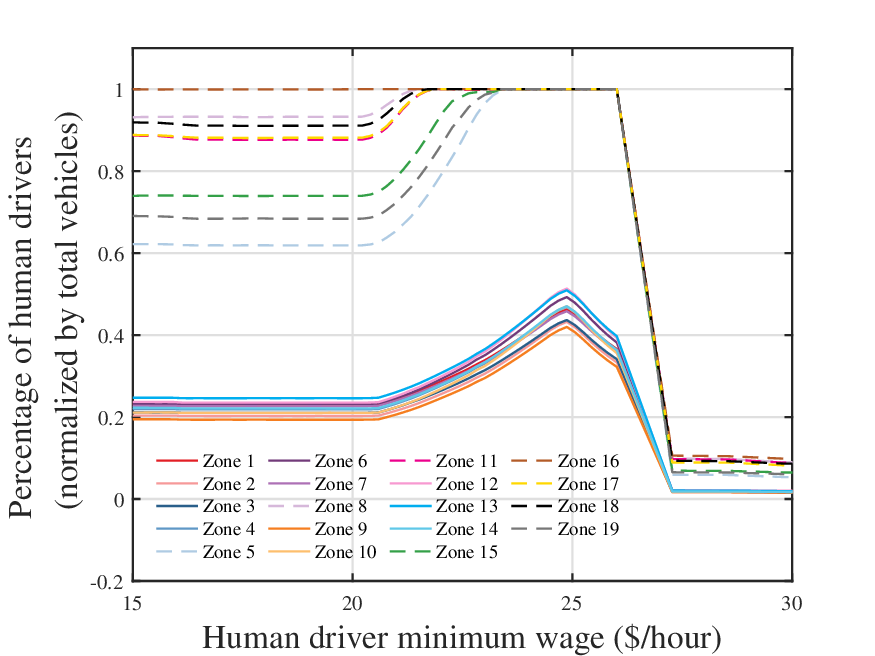}
    \caption{Percentage of human driver in each zone under varying minimum driver payment $\it{q_{min}}$ (normalized by total vehicles)}
    \label{allidleHDV_r2}
    \end{minipage}
    \begin{minipage}[!htbp]{0.49\textwidth}
    \includegraphics[height=6.4cm,width=8.7cm]{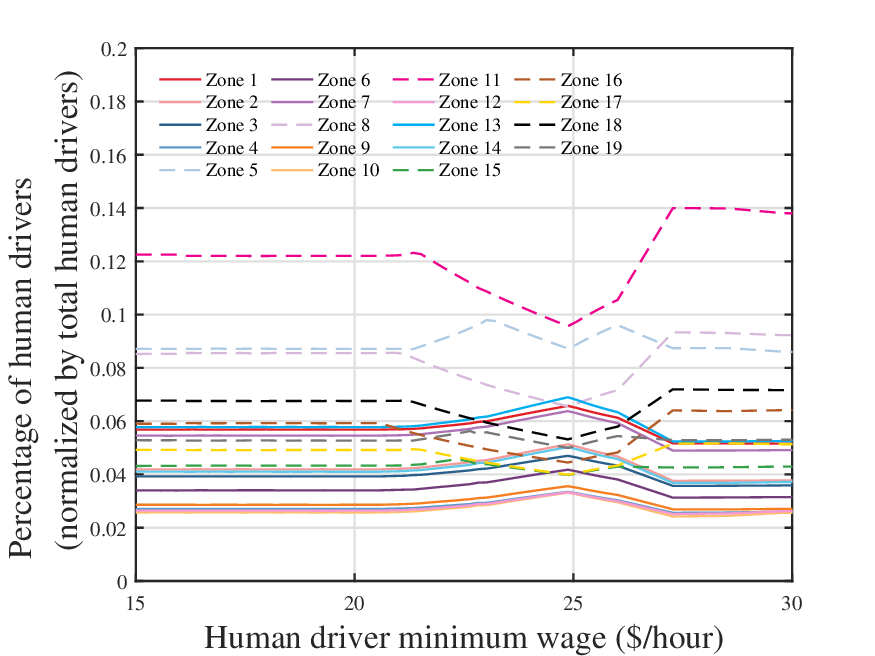}
    \caption{Percentage of human driver in each zone under varying minimum driver payment $\it{q_{min}}$ (normalized by total human drivers)}
    \label{allidleHDV_r22}
    \end{minipage}
    \end{figure}


In the rest of this subsection, we fix AV cost at $C=\$\,26$/hour and vary the minimum driver payment $q_{min}$ to investigate how the minimum wage regulation affects the TNC market. The profit-maximization problem (\ref{optimalpricing_trip_regulated}) can be solved through dual decomposition in a similar way as proposed in Algorithm \ref{algorithm1}. The algorithm can provide an approximate solution which satisfies all the constraints in (\ref{profit_constraints_regulated}). The performance of the approximate solution is shown in Figure \ref{Comparison_relax_r2}, where the maximum gap between (\ref{optimalpricing_trip_regulated}) and its relaxed problem is 5.21\,\%. The detailed results of the simulation study are summarized in Figure \ref{Comparison_relax_r2}-\ref{allidleHDV_r22}. The key insights are summarized below:
\begin{enumerate}
    \item In certain regime, i.e., $\$\,21$/hour $\leq q_{min}<\,\$\,24.8$/hour, a higher minimum wage for human driver will motivate the platform to hire {\em more} human drivers and {\em fewer} AVs (Figure \ref{totalnovehicle_r2}). This seems counter-intuitive as it indicates that the platform recruits more human drivers even if each human driver is more expensive. We comment that this is because in the labor market for human drivers, the platform is a monopsony that will hire fewer employees compared to firms in the competitive market setting in order to maximize its profit \cite{li2019regulating}. A minimum wage on human drivers will curb the market power of the TNC platform and force the platform to hire more human drivers so that the market outcome is closer to that of the competitive equilibrium. This is in line with our previous work \cite{li2019regulating}, which provides a rigorous proof for this conclusion and showed that it always holds under arbitrary model parameters with an aggregate equilibrium model. We conjecture (without a rigorous proof) that the same conclusion should also apply here and the above observation is robust under a large range of model parameters.
    
    \item Beyond certain point, i.e., $q_{min}\geq\,\$\,24.8$/hour, the increasing trend of human-driver supply is reversed, and the increasing minimum wage triggers a paradigm shift where the platform sharply reduces the number of human drivers and replaces them with AVs, until the majority of the human drivers are removed from the TNC fleet. This is an intuitive result: when the minimum wage is too high, it increases the cost of human labor for the platform, and it prefers to use AVs instead of human drivers. It indicates that minimum wage has to be carefully selected in the mixed environment to avoid massive loss of job opportunities for human drivers, which is more likely to occur in the mixed environment because of the availability of more affordable supply (e..g, AVs).
    \item Similar to the results in Section \ref{subsection5_1}, the platform prefers to concentrate AVs in high-demand areas and dispatch human drivers to suburban areas under the minimum wage. This is evidenced by Figure \ref{idleAV_r2}-\ref{allidleHDV_r22}, where more human drivers are placed in suburban areas than that in urban areas. However, different from Section \ref{subsection5_1}, even if the increase of minimum wage leads to a higher concentration of human drivers in suburban areas, the earnings of each driver increases. This reveals that drivers are effectively protected by the minimum wage when $\$\,21$/hour $\leq q_{min}<\,\$\,24.8$/hour.
\end{enumerate}

To summarize, a minimum wage on human drivers will effectively protect for-hire human drivers from the negative impacts of AV deployment under certain regime. However, an overly protective minimum wage will trigger a paradigm shift where platform sharply reduces the number of human drivers and replaces most of them with AVs. Therefore, minimum wage should be exercised with discretion in the mixed environment to avoid massive loss of job opportunities for human drivers.

\subsection{Impacts of Pickup Restriction on AVs}
\label{sectin 5.3}
To improve the transport equity of ride-hailing services for passengers, we consider a restrictive pickup policy that prevents AVs from picking up passengers in urban areas, while allowing them to drop off passengers in these areas. This is consistent with existing regulations in many cities: in New York City, green cabs are not permitted to street-hail passengers in South Manhattan; in Hong Kong, green taxi is only available for hire in New Territory, and blue taxi is only available in Lantau Island, where as red taxi can pick up passengers everywhere in the city. Admittedly, this will increase the deadhead miles because taxis that drop off a passenger in restricted areas will have to travel back to permitted areas with empty seats. However, it will motivate drivers to redistribute to remote areas, which will simultaneously improve traffic condition in urban core and enhance service quality for passengers in underserved communities.

In the case study, suburban areas include zone 5, zone 8, zone 11, zone 15, zone 16, zone 17, zone 18, and zone 19, and the proposed policy enforces that AVs can only pick up passengers from these zones. This can be easily incorporated into (\ref{optimalpricing_trip}) by setting idle AVs outside of these areas to be 0, thus passengers starting from urban zones can only hail a ride with human drivers. Under the proposed policy, the profit-maximizing problem of the platform is similar to (\ref{optimalpricing_trip}), and we can still solve it through Algorithm \ref{algorithm1}.

\begin{figure}[h!]
    \centering
    \begin{minipage}[h]{0.48\textwidth} 
    \includegraphics[width=1\textwidth]{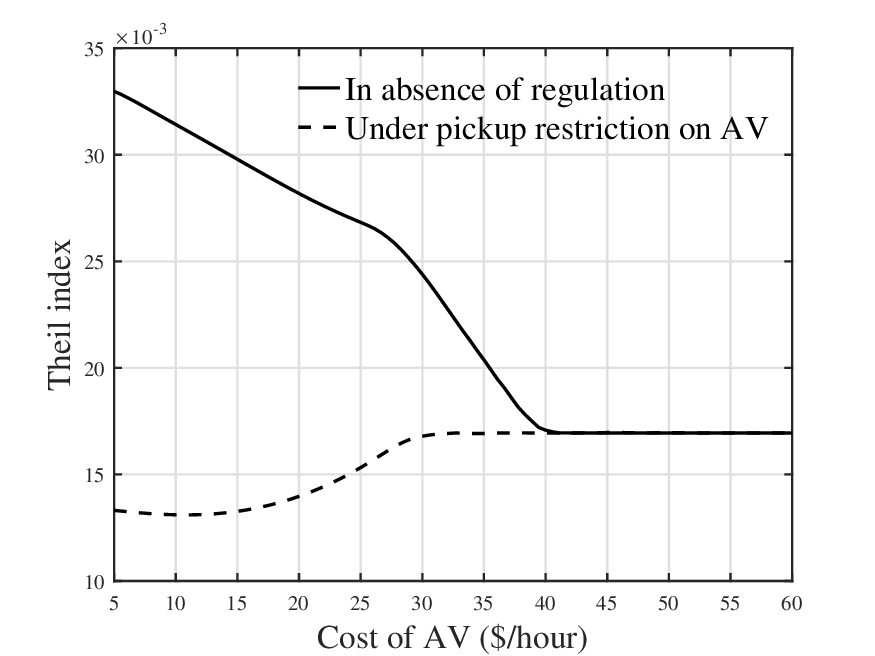}
    \caption{Theil coefficient for expected travel cost with and without pickup restriction on AV}
    \label{cvcomr3}
    \end{minipage}
    \begin{minipage}[h]{0.48\textwidth} 
    \includegraphics[width=1\textwidth]{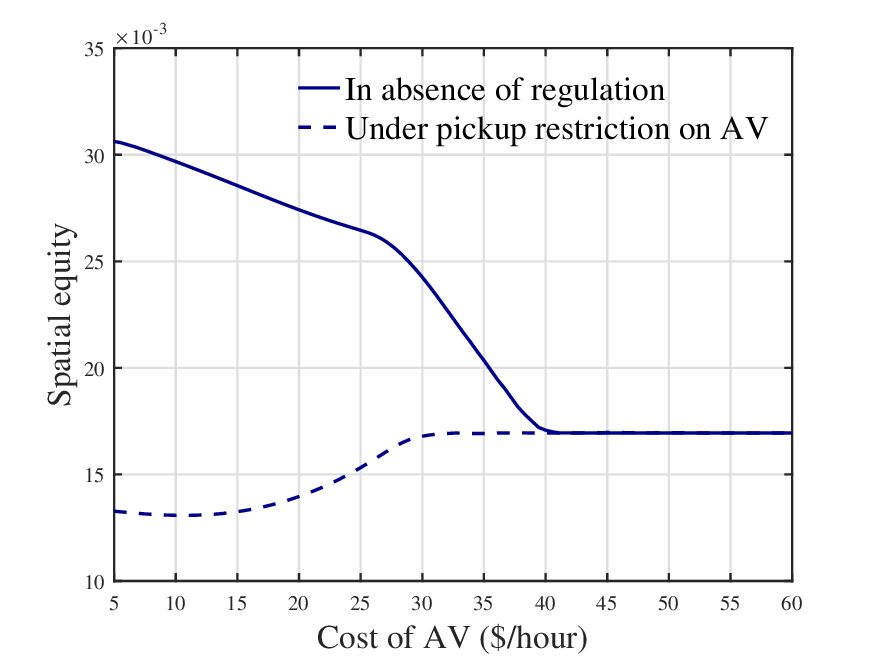}
    \caption{Spatial equity for expected travel cost with and without pickup restriction on AV}
    \label{Spatial_equity}
    \end{minipage}
    \begin{minipage}[h]{0.48\textwidth} 
    \includegraphics[width=1\textwidth]{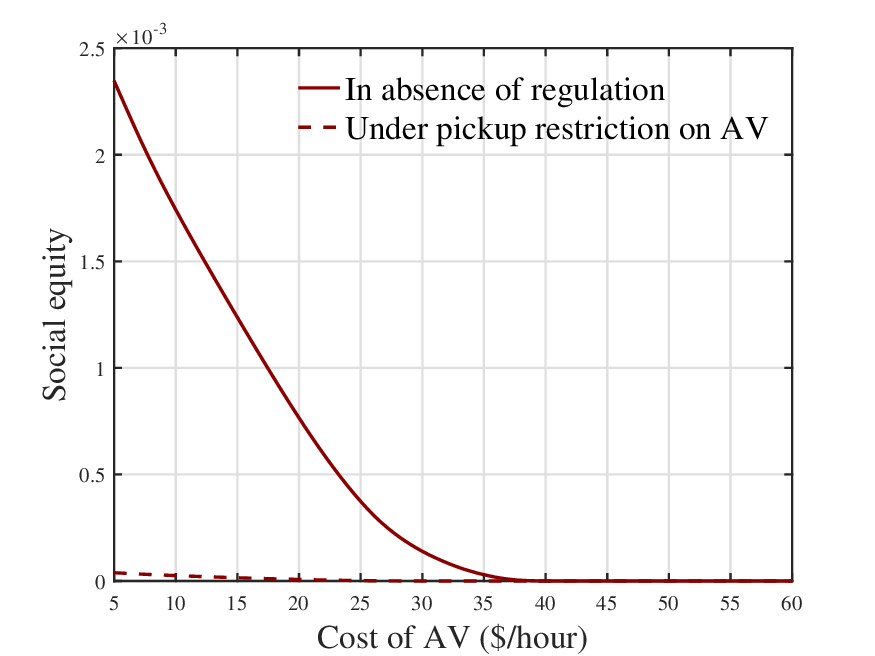}
    \caption{Social equity for expected travel cost with and without pickup restriction on AV}
    \label{Vertical_equity}
    \end{minipage}
    \begin{minipage}[h]{0.49\textwidth}
    \centering
    \includegraphics[width=1\textwidth]{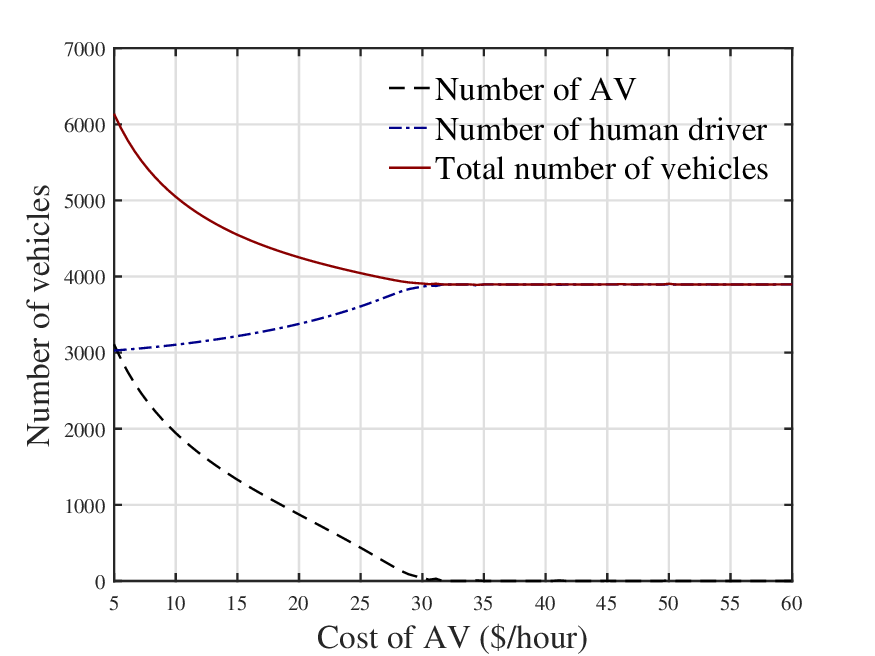}
    \caption{Total number of vehicles in the TNC fleet under pickup restriction of AV}
    \label{aggtotalvehicler3}   
    \end{minipage}
    \begin{minipage}[h]{0.49\textwidth} 
    \includegraphics[width=1\textwidth]{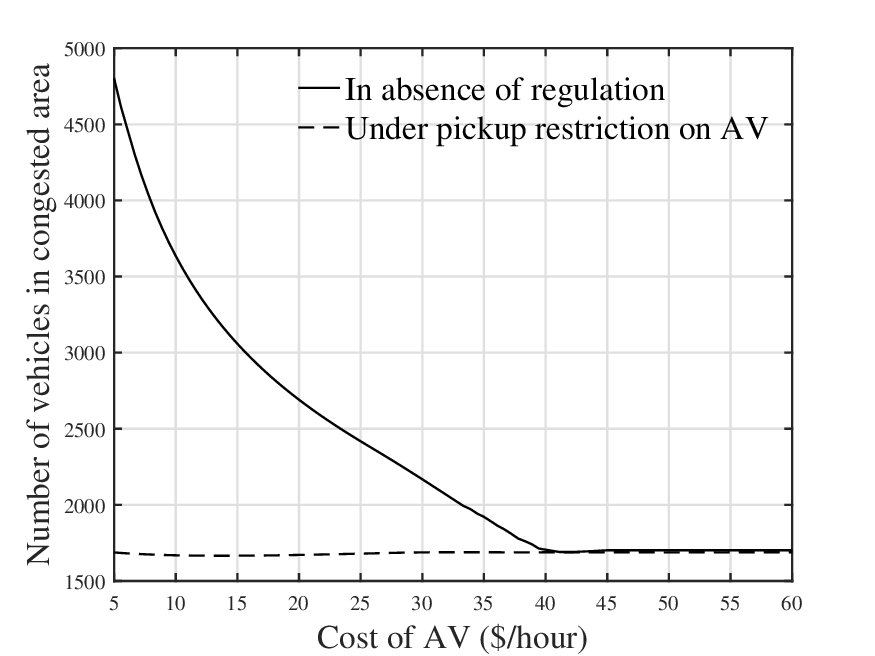}
    \caption{Number of vehicles in congested area}
    \label{recordnc}
    \end{minipage}
    \begin{minipage}[h]{0.49\textwidth} 
    \includegraphics[width=1\textwidth]{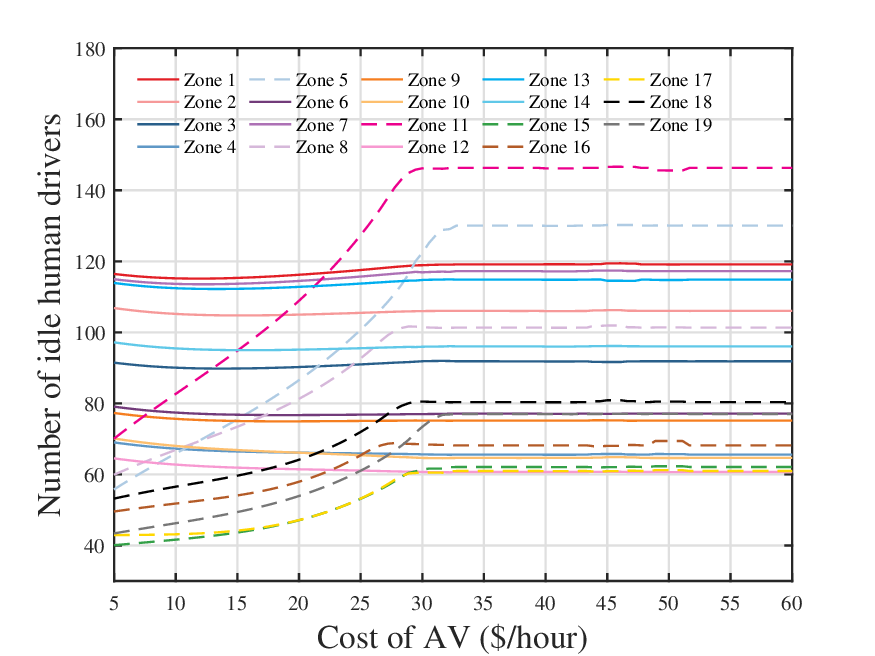}
    \caption{Number of idle human drivers in each zone under pickup restriction of AV}
    \label{idleHDVr3}
    \end{minipage}
\end{figure}

\begin{figure}[h!]
    \centering
    \begin{minipage}[h]{0.49\textwidth} 
    \centering
    \includegraphics[width=1\textwidth]{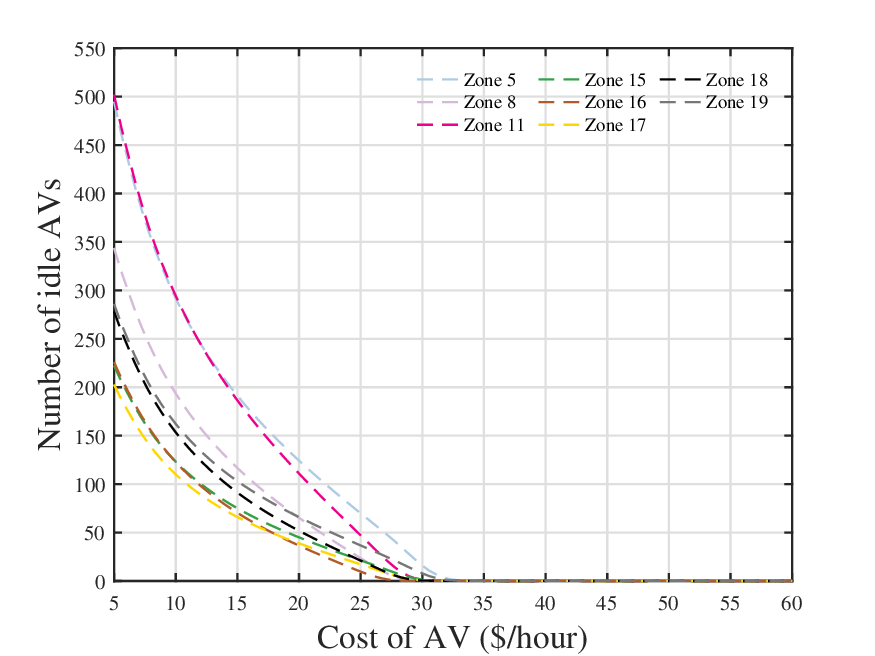}
    \caption{Number of idle AVs in each zone under pickup restriction of AV}
    \label{idleavr3}
    \end{minipage}
    \begin{minipage}[h]{0.49\textwidth} 
    \includegraphics[width=1\textwidth]{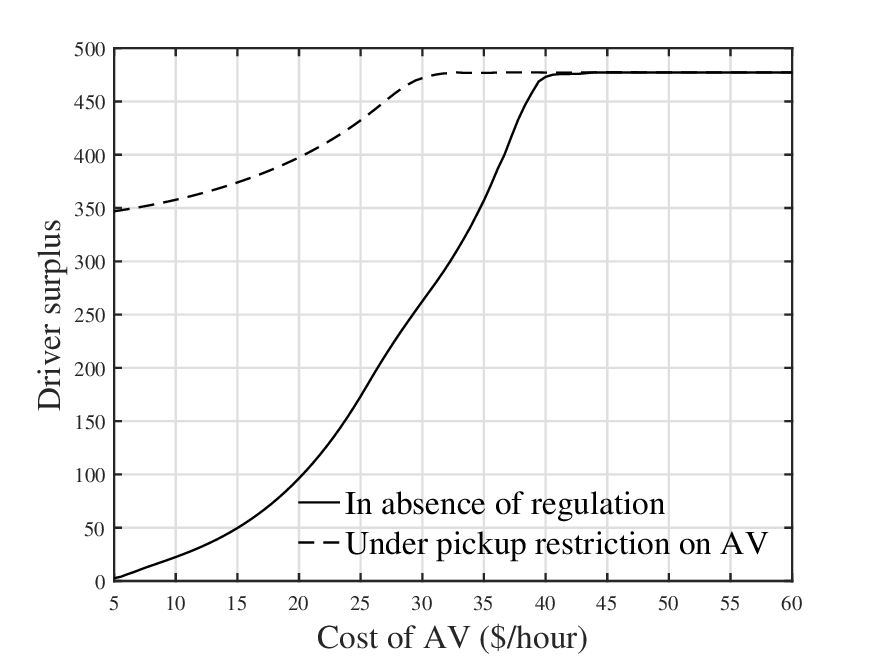}
    \caption{Comparison of driver surplus with and without pickup restriction of AV}
    \label{driversurplusr3}
    \end{minipage}
    \begin{minipage}[h]{0.49\textwidth} 
    \includegraphics[height=6.2cm,width=8.5cm]{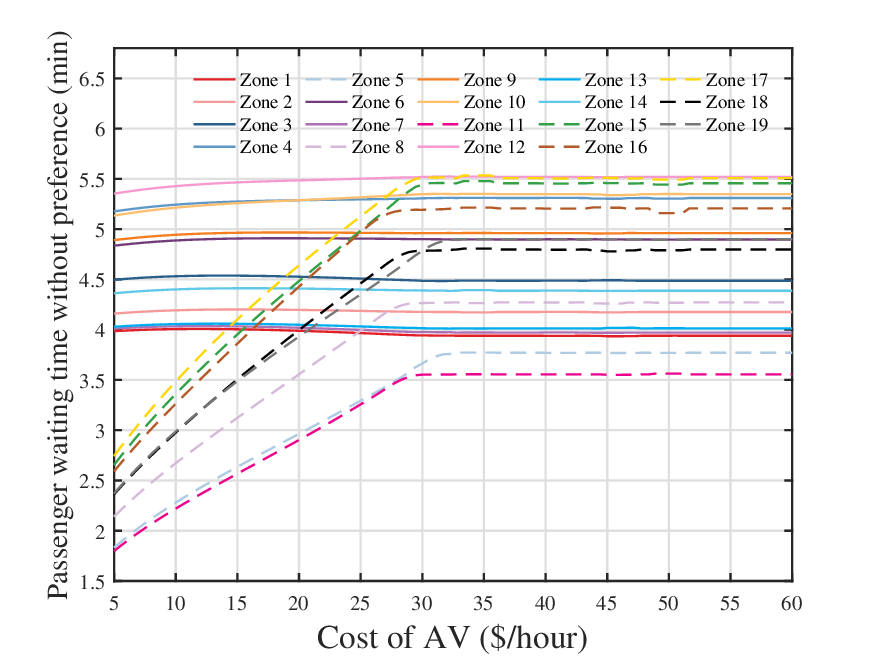}
    \caption{Passenger without preference waiting time in each zone under pickup restriction of AV}
    \label{passengerwaitr3}
    \end{minipage}
    \begin{minipage}[h]{0.49\textwidth} 
    \includegraphics[height=6.2cm,width=8.5cm]{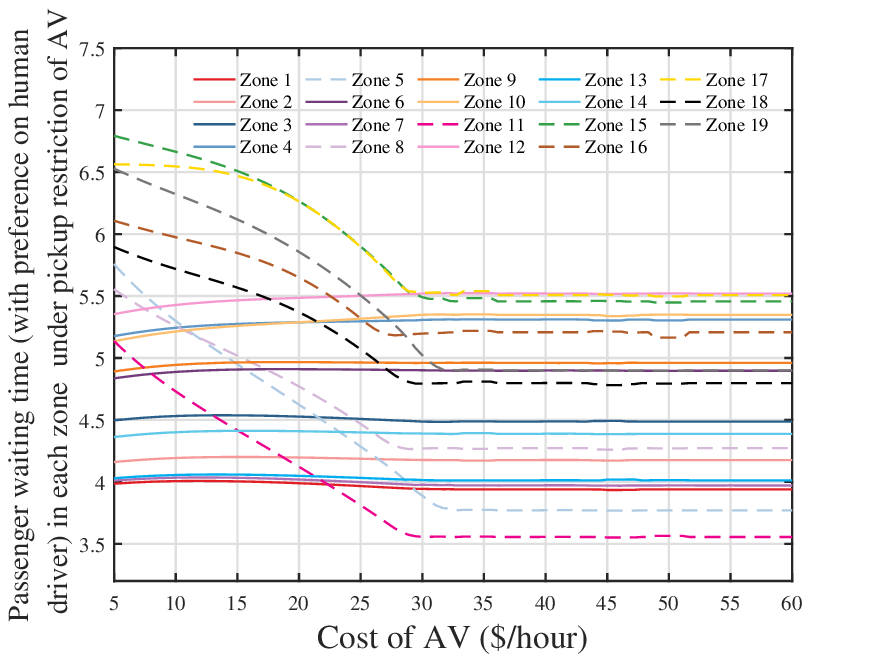}
    \caption{Passenger with preference waiting time in each zone under pickup restriction of AV}
    \label{passengerwaitr4}
    \end{minipage}
    \begin{minipage}[h]{0.49\textwidth} 
    \includegraphics[width=1\textwidth]{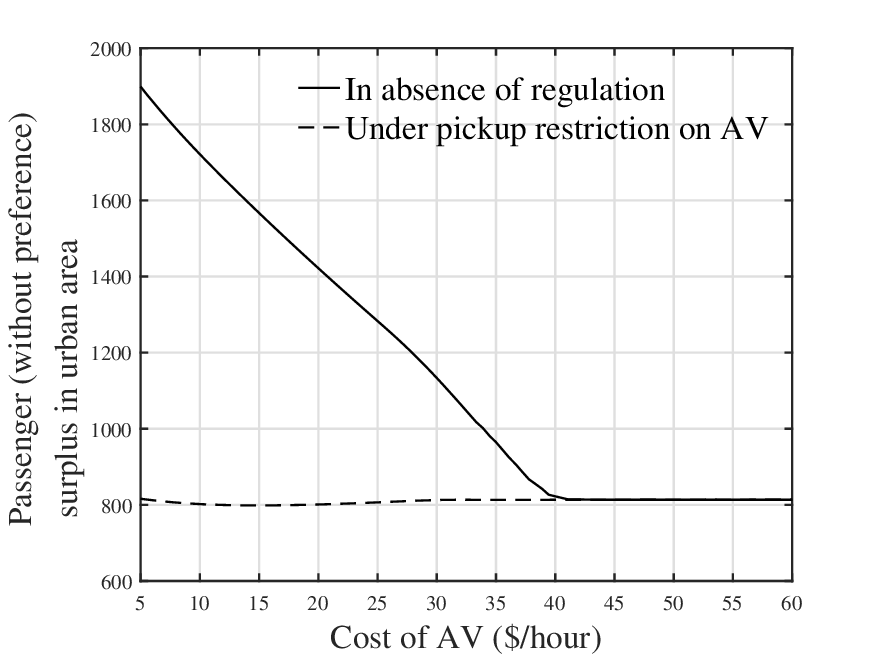}
    \caption{Comparison of passenger (without preference) surplus in urban area with and without pickup restriction of AV}
    \label{passengersurpluscomurbanwithout}
    \end{minipage}
    \begin{minipage}[h]{0.49\textwidth} 
    \includegraphics[width=1\textwidth]{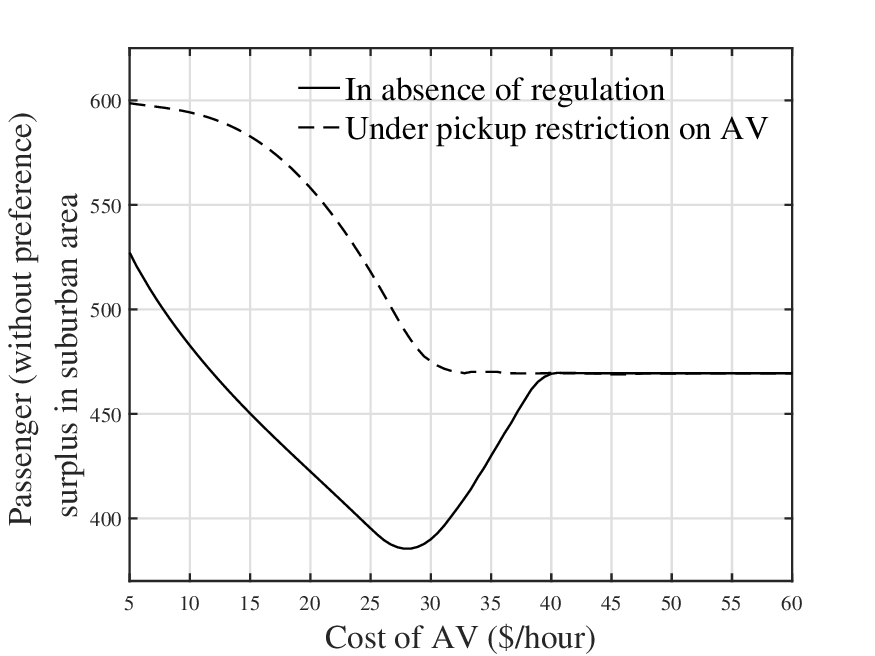}
    \caption{Comparison of passenger (without preference) surplus in suburban area with and without pickup restriction of AV}
    \label{passengersurpluscomsuburbanwithout}
    \end{minipage}
\end{figure}

\begin{figure}[t]
    \centering
    \begin{minipage}[h]{0.32\textwidth} 
    \includegraphics[width=0.9\textwidth]{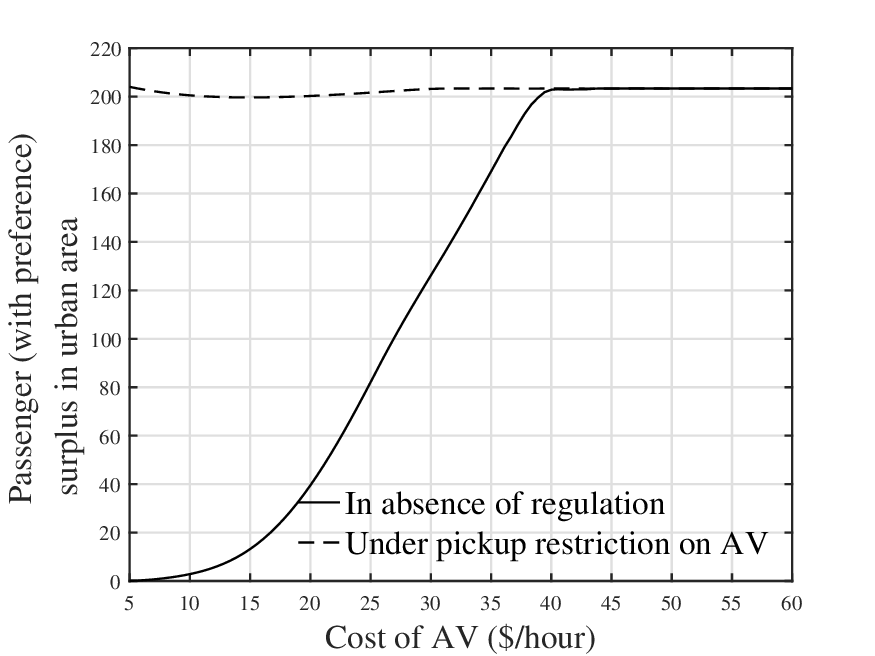}
    \caption{Comparison of passenger (with preference) surplus in urban area with and without pickup restriction of AV}
    \label{passengersurpluscomurbanwith}
    \end{minipage}
    \begin{minipage}[h]{0.32\textwidth} 
    \includegraphics[width=0.9\textwidth]{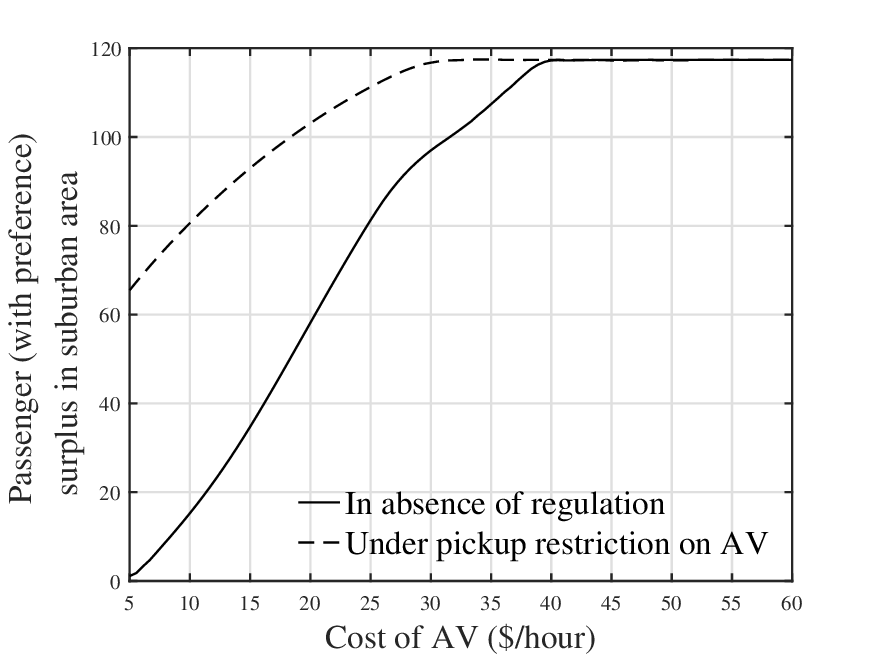}
    \caption{Comparison of passenger (with preference) surplus in suburban area with and without pickup restriction of AV}
    \label{passengersurpluscomsuburbanwith}
    \end{minipage}
    \begin{minipage}[h]{0.32\textwidth}
    \centering
    \includegraphics[width=0.9\textwidth]{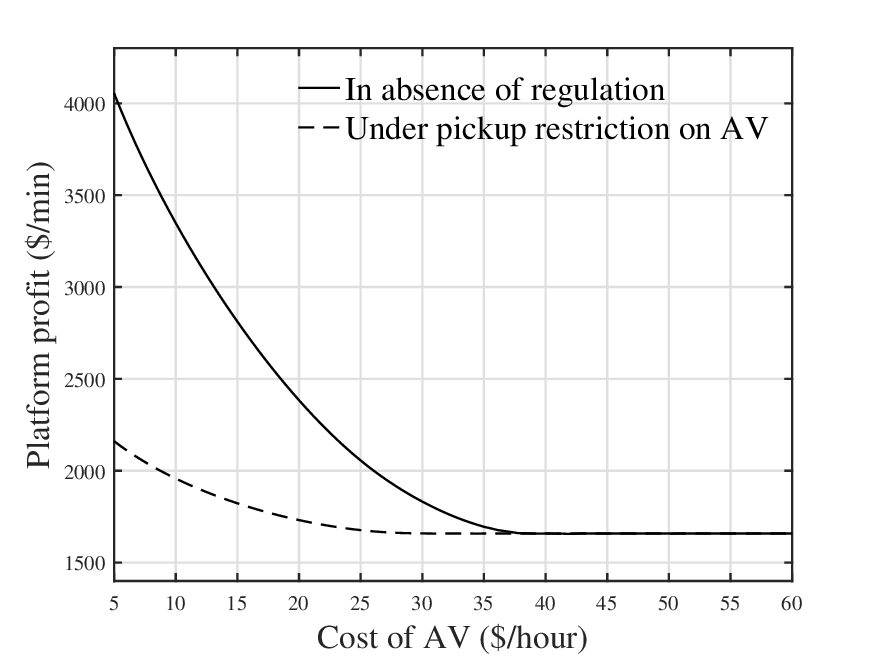}
    \caption{Comparison of TNC platform profit with and without pickup restriction of AV}
    \label{profitcomr3}
    \end{minipage}
\end{figure}

Figures \ref{cvcomr3}-\ref{aggtotalvehicler3} present the outcomes of the TNC market under the proposed AV policy when the cost of AV changes. We summarize the key insights as follows:

\begin{itemize}
   \item The proposed AV policy can effectively improve transport equity for passengers from both spatial and social perspective. Under the regulation policy, the overall Theil index has  significantly reduced compared to the unregulated case (Figure \ref{cvcomr3}). At the same time, the spatial inequity gap reduces when the AV cost reduces (Figure \ref{Spatial_equity}), whereas the social inequality gap slightly increases as AV cost reduces (Figure \ref{Vertical_equity}), but in both cases, when compared with the unregulated case, the corresponding index at the same AV cost is substantially smaller. This indicates that under the proposed regulation policy, platforms will prioritize AV deployment in underserved areas, which positively affects spatial transport equity. In the meanwhile, the platform will hire more human drivers in the overall market, which benefits disadvantaged groups who has special needs and therefore have to take a ride-hailing trip served by human drivers.
   
   \item Interestingly, the proposed policy can limit the total number of ride-hailing vehicles compared to the unregulated case and alleviate the traffic congestion associated with AV adoption. It can be observed that the increase of total fleet size (Figure \ref{aggtotalvehicler3}) as AV cost decreases under the pickup restriction  policy is much smaller as compared to that in the unregulated scenario (Figure \ref{totalnovehicle}). Furthermore, it is obvious that the total number of vehicles in the congested areas of the city is also much smaller in the regulated case as compared to the unregulated case (Figure \ref{recordnc}). This indicates that the aforementioned regulatory policy has the potential to restrict the increase of vehicle fleet and mitigate the traffic congestion in core urban areas.

   \item The proposed AV policy induces significant benefits for human drivers. As AVs are restricted to pick up passengers in suburban areas, human drivers becomes indispensable for the TNC platform in urban areas. Therefore, as AV cost reduces, AV deployment scales up, but they are primarily used to replace human drivers in suburbs (see Figure \ref{idleHDVr3} and \ref{idleavr3}). As a consequence, many human drivers will retain their jobs even if AV cost is low (Figure \ref{aggtotalvehicler3}), and surplus of human drivers is much higher compared to the unregulated case (Figure \ref{driversurplusr3}).
   
   \item The impacts of AV policies on passengers are are nuanced and may depend on where passengers start their trips and what class of passenger we refer to. Under the proposed AV policy, reduced AV cost will not significant influence the passengers in urban areas regardless of passenger classes (see Figure \ref{passengerwaitr3} and \ref{passengerwaitr4}). On the other hand, the decreased AV cost will lead to lower waiting time for passengers in remote areas who has no preferences (see Figure \ref{passengerwaitr3}), but higher waiting time for passengers in remote areas who has preference over human drivers (see Figure \ref{passengerwaitr4}). This is because the proliferation of AVs in remote areas will enhance the competition between AVs and human drivers, incentivizing more human drivers to move to the urban area. Thus, passengers living in remote areas and preferring human drivers will find it more challenging to hail a ride. On the other hand, when comparing the case with and without regulation under the same AV cost, it is clear that passengers with a preference for human drivers all benefit from the regulation, regardless of whether they are in urban or remote areas (See Figure \ref{passengersurpluscomurbanwith} and \ref{passengersurpluscomsuburbanwith}). Simultaneously, for passengers without a preference, those in remote areas benefit from the regulation (see Figure \ref{passengersurpluscomsuburbanwithout}), but those in urban areas may be negatively affected (see Figure \ref{passengersurpluscomurbanwithout}). This is quite intuitive: the pickup restriction limits the use of AVs in urban areas, which could have been available for this class of users if the regulation policy were not in place.
   
\end{itemize}

In summary, by restricting the pickup locations of AVs, the platform is motivated to prioritize the deployment of AVs in suburban areas. This not only improves the transport equity over the transportation network for passengers, but simultaneously protect human drivers from being completely replaced by AVs. However, the proposed AV policy may negatively affect passengers in the urban areas and lead to profit reduction for the TNC platform (Figure \ref{profitcomr3}). This reflects a trade-off between different policy objectives (e.g., equity v.s. service quality), which warrants further studies in future works.

\begin{remark}
The restrictive pickup policy considered in this paper may impede the popularization of AVs in the ride-hailing industry. This is evident as the number of AVs in Figure \ref{aggtotalvehicler3} (representing the regulated case) is smaller than that in Figure \ref{totalnovehicle} (representing the unregulated case). This prompts the question: should the government still adopt such a policy in practice, despite its potential to hinder the uptake of AVs? On this matter, we clarify that the proliferation of AVs on ride-hailing platforms carries both benefits and drawbacks. The benefits are evident: introducing AVs into the ride-hailing business can reduce operational costs compared to human drivers, enhance controllability in vehicle repositioning, improve service responsiveness, and potentially lower passenger fares. However, the rapid expansion of AVs on ride-hailing platforms also introduces significant concerns, including job losses for human drivers, inequitable spatial distribution of vehicles, increased traffic congestion, and uneven distribution of benefits among various passenger types. It remains unclear whether the benefits of AVs will necessarily surpass their negative externalities. Therefore, any policy decision should be made after a comprehensive evaluation of its overall positive and negative impacts.  Our model offers a systematic framework for assessing these benefits and negative externalities, thereby providing a valuable tool for governments to make informed policy decisions.
\end{remark}

\section{Conclusion}
This paper studies a TNC platform with a mixed fleet of AVs and human drivers over the transportation network and investigates how autonomous vehicles will affect human drivers and passengers in a regulated TNC market. A spatial equilibrium model is proposed to captures the spatial distribution of AVs, the flow of occupied and idle vehicles, the flow of passengers, and traffic congestion over the transportation network. The overall problem is cast as a non-concave problem, where the TNC platform determines the spatial prices, the driver wage, and the placement of idle AVs to maximize its profit. We show that the solution to the network equilibrium constraints exits under mild conditions, and we develop a solution algorithm to calculate an approximate solution to the non-concave problem with a verifiable performance guarantee.

Based on the proposed model, we first investigate how the cost of AVs will affect the market equilibrium. We show that as the AV cost reduces, the platform hires more AVs and fewer human drivers, which reduces the earning and job opportunities of human drivers. We also observe that the TNC platform prioritizes AV deployment in the urban core as there is higher passenger demand. This renders competition between AVs and human drivers in these areas and pushes human drivers to relocate to suburbs, which further reduces their earning opportunities and intensifies the spatial inequity of service quality among passengers. To mitigate these concerns, we conduct further policy analysis in two directions: (a) to protect human drivers, we investigate the impacts of a wage floor on human driver earnings. We show that a carefully selected wage floor may benefit human drivers without significantly affecting passengers. However, we find that there exists a threshold beyond which the minimum wage will trigger a paradigm shift and the platform will completely replace all human drivers with AVs; (b) to improve the spatial equity for passengers, we considered a pickup restriction policy, which prohibits AVs from picking up passenger in urban areas. Our study shows that the transport inequity can be improved and the fast increase of total fleet size due to the proliferation of AVs can be restricted with the implementation of the pickup restriction policy, while it also cut down the surplus of passengers in urban areas and decreases the profit of the TNC platform. These insights provide guidance for the policymaker to understand the potential impact of AVs and make preparations for relevant regulatory interventions.

In the future, we will extend the proposed framework to characterize the temporal dynamics of passenger and driver, investigate the impacts of transport electrification, examine the interaction between TNC and other transport models, and develop other regulatory policies that balance the trade-off between service quality and equity of the TNC market.

\section*{Acknowledgments}  {This research was supported by National Science Foundation of China under project 72201225 and Hong Kong Research Grants Council under project 16202922.}

\bibliographystyle{unsrt}
\bibliography{ref} 
\appendix

\section{Correspondences Between Postal Code and Zone Number}
\label{Appendix B}
The following Table \ref{mytable} describes the correspondences between postal code and zone number:
    \begin{table}[!h]
    \setlength{\abovecaptionskip}{0.3cm}
    \setlength{\belowcaptionskip}{0.5cm}
        \centering
        \caption{The correspondence between zone number and postal code}
        \begin{tabular}{c|c}
            1 &\quad94104\tablefootnote{We integrate the zones with postal codes 94104, 94105, 94111 as an aggregated zone with postal code 94104.} \\
            2 & 94103 \\
            3 & 94109 \\
            4 & 94115 \\
            5 & 94118 \\
            6 & 94123 \\
            7 &\quad94108\tablefootnote{We integrate the zones with postal codes 94108 and 94133 as an aggregated zone with postal code 94108.} \\
            8 & 94121 \\
            9 & 94102 \\
            10 & 94117 \\
            11 & 94122 \\
            12 & 94114 \\
            13 & 94107 \\
            14 & 94110 \\
            15 & 94131 \\
            16 & 94116 \\
            17 & 94124 \\
            18 & 94132 \\
            19 & 94112
        \end{tabular}
    \label{mytable}
    \end{table}
    
    \newpage

\end{document}

%% file: format.tex
\usepackage{epsfig}
\usepackage{subfigure}
\usepackage{amssymb,amsmath,amsfonts,layout,graphicx}
\usepackage{makeidx}
\usepackage{babel}
\usepackage{tikz}
\usepackage{sublabel}
\usepackage{cases}
\usepackage{algorithm}
\usepackage{algorithmic}

\usepackage
[
        letterpaper,
        left=1.91cm,
        right=1.91cm,
        top=1.91cm,
        bottom=2cm,
]
{geometry}

\newtheorem{lemma}{Lemma}

\newtheorem{proposition}{Proposition}

\newtheorem{remark}{Remark}

\newcommand{\x}{{\bf x}}